\documentclass[11pt]{article}
\usepackage{amsfonts,latexsym,amsmath,amscd}
\newcommand{\beq}{\begin{equation}}
\newcommand{\eeq}{\end{equation}}
\newcommand{\beqa}{\begin{eqnarray}}
\newcommand{\eeqa}{\end{eqnarray}}
\newcommand{\beaa}{\begin{eqnarray*}}
\newcommand{\ben}{\begin{eqnarray*}}
\newcommand{\eaa}{\end{eqnarray*}}
\newcommand{\een}{\end{eqnarray*}}

\newcommand{\bal}{\begin{align}}
\newcommand{\eal}{\end{align}}

\newcommand \nc {\newcommand}
\nc \proof {\noindent {\em{Proof.\/ }}}
\nc \qed {\hfill $\Box$}
\newtheorem{theorem}{Theorem}[section]
\newtheorem{lemma}[theorem]{Lemma}
\newtheorem{proposition}[theorem]{Proposition}
\newtheorem{corollary}[theorem]{Corollary}
\newtheorem{definition}[theorem]{Definition}
\newtheorem{example}[theorem]{Example}
\newtheorem{remark}[theorem]{Remark}
\newtheorem{conjecture}[theorem]{Conjecture}
\newtheorem{question}[theorem]{Question}
\nc \bth[1] {\begin{theorem}\label{t#1} }
\nc \ble[1] {\begin{lemma}\label{l#1} }
\nc \bpr[1] {\begin{proposition}\label{p#1} }
\nc \bco[1] {\begin{corollary}\label{c#1} }
\nc \bde[1] {\begin{definition}\label{d#1}\rm }
\nc \bex[1] {\begin{example}\label{e#1}\rm }
\nc \bre[1] {\begin{remark}\label{r#1}\rm }
\nc \bcon[1] {\begin{conjecture}\label{con#1}\rm }
\nc \bque[1] {\begin{question}\label{que#1}\rm }
\nc {\eth} {\end{theorem}}
\nc {\ele} {\end{lemma}}
\nc {\epr} {\end{proposition}}
\nc {\eco} {\end{corollary}}
\nc {\ede} {\end{definition}}
\nc {\eex} {\end{example}}
\nc {\ere} {\end{remark}}
\nc {\econ} {\end{conjecture}}
\nc {\eque} {\end{question}}
\nc \thref[1]{Theorem \ref{t#1}}
\nc \leref[1]{Lemma \ref{l#1}}
\nc \prref[1]{Proposition \ref{p#1}}
\nc \coref[1]{Corollary \ref{c#1}}
\nc \deref[1]{Definition \ref{d#1}}
\nc \exref[1]{Example \ref{e#1}}
\nc \reref[1]{Remark \ref{r#1}}

\def \( {{\left(}}
\def \) {{\right)}}

\def \[  {{\left[}}
\def \]  {{\right]}}

\def\tensor{\otimes}

\def\x{\xi}

\begin{document}
\title{Regularity of Dirac-harmonic maps}
\author{Changyou Wang\\
Department of Mathematics, University of Kentucky\\
Lexington, KY 40506, USA\\
{\it cywang@ms.uky.edu}\\
\\
Deliang Xu\\Department of Mathematics, Shanghai Jiaotong University\\
Shanghai 200240, P. R. China\\
{\it dlxu@sjtu.edu.cn}}
\date{}
\maketitle

\begin{abstract}

For any $n$-dimensional compact spin Riemannian manifold $M$ with a given spin structure and a spinor bundle $\Sigma M$,
and any compact Riemannian manifold $N$, we show an $\epsilon$-regularity
theorem for weakly Dirac-harmonic maps $(\phi,\psi):M\tensor\Sigma M
\to N\tensor\phi^*TN$. As a consequence, any weakly Dirac-harmonic map is proven to be smooth when $n=2$.
A weak convergence theorem for approximate Dirac-harmonic maps is established when $n=2$.
For $n\ge 3$, we introduce the notation of
stationary Dirac-harmonic maps and obtain a Liouville theorem for stationary Dirac-harmonic maps
in $\mathbb R^n$. If, additions, $\psi\in W^{1,p}$ for some
$p>\frac{2n}{3}$, then we obtain an energy monotonicity formula and prove a partial regularity
theorem for any such a stationary Dirac-harmonic map.
\end{abstract}


\section{Introduction}

The notation of Dirac-harmonic maps is inspired by the supersymmetric nonlinear sigma
model from the quantum field theory \cite{deligne}, and is a very natural and interesting
extension of harmonic maps. In a series of papers \cite{chen-jost-li-wang1,
chen-jost-li-wang2}, Chen-Jost-Li-Wang recently introduced the subject of Dirac-harmonic
maps and studied some analytic aspects of Dirac-harmonic maps from a spin Riemann surface into
another Riemannian manifold.
In order to review some of the main theorems of \cite{chen-jost-li-wang1, chen-jost-li-wang2} and
motive the aim of this paper, let's briefly describe the mathematical framework given by
\cite{chen-jost-li-wang1, chen-jost-li-wang2}.

For $n\ge 2$, let $\left(M,g\right)$ be a compact $n$-dimensional spin Riemannian
manifold with a given spin structure and an associated spinor bundle $\Sigma(=\Sigma M)$,
and $\left( N,h\right)$ be a compact $k$-dimensional Riemannian manifold without
boundary. By Nash's theorem, we may assume that $\left(N,h\right)$ is isometrically
embedded into an Euclidean space $\mathbb{R}^K$ for a sufficiently large $K$.
Let $\nabla^M$ and $\nabla^N$ be the Levi-Civita connection on $M$ and $N$ respectively.
Let $\langle\cdot,\cdot\rangle$ be a Hermitian metric on $\Sigma$ (a
complex vector bundle of complex dimension $n$), and $\nabla^{\Sigma}$ be
the Levi-Civita connection on $\Sigma$ compatible with the metrics $\langle\cdot,\cdot\rangle$
and $g$.  For a map $\phi:M\to N$, let $\phi^{*}TN$ denote the pull-back bundle of $TN$ by $\phi$ that is
equipped with the pull-back metric $\phi^*h$ and the connection $\nabla^{\phi^*TN}$.
On the twist bundle $\Sigma\tensor \phi^*TN$, there is a metric, still denoted as $\langle\cdot,\cdot\rangle$,
induced from the $\langle\cdot,\cdot\rangle$ and $\phi^*h$.
There is also a Levi-Civita connection $\overline\nabla$ on $\Sigma\tensor\phi^*TN$
induced from $\nabla^\Sigma$ and $\nabla^{\phi^{*}TN}$.  The Dirac
operator $\mathcal D$ along $\phi$ is defined as follows. For any section
$\psi\in\Gamma(\Sigma\tensor\phi^*TN)$,
\beq{}\label{dirac-op}
\psi^i=f_\alpha\circ\overline\nabla_{f_\alpha}\psi,
\eeq
where $\{f_\alpha\}_{\alpha=1}^n$ is a local orthonormal frame on $M$, and
$\circ: TM\tensor_{\mathbb C}\Sigma\to\Sigma$ is the Clifford multiplication.
More precisely, if we write $\psi$ in the local coordinate as
$\psi=\psi^i\tensor\frac{\partial}{\partial y_i}(\phi)$,
where $\psi^i\in\Gamma\Sigma$ is a section of $\Sigma$ for $1\le i\le k$ and $\{\frac{\partial}{\partial y_i}\}_{i=1}^k$ is
a local coordinate frame on $N$, then
\beq{}\label{dirac-op1}
D\hspace{-0.25cm}/\hspace{0.1cm}\psi^i=\partial\hspace{-0.2cm}/\hspace{0.1cm} \psi^i\tensor \frac{\partial}{\partial y_i}(\phi)
+\left(f_\alpha\circ\psi^i\right)\tensor\nabla_{f_\alpha}^{\phi^*TN}\left(\frac{\partial}{\partial y_i}(\phi)\right),
\eeq
where $\partial\hspace{-0.2cm}/\hspace{0.1cm} =f_\alpha\circ\nabla^\Sigma_{f_\alpha}$ is the standard Dirac operator on the spin
bundle $\Sigma$.

The Dirac-harmonic energy functional was first
introduced by Chen-Jost-Li-Wang in \cite{chen-jost-li-wang1, chen-jost-li-wang2}.
\beq{}\label{dirac-harmonic-eng}
L\left(\phi,\psi\right) =\int_M\left[\left\vert d \phi\right\vert^2+\left\langle
\psi,D\hspace{-0.25cm}/\hspace{0.1cm}\psi\right\rangle\right]dv_g
=\int_M \left[g^{\alpha\beta}h_{ij}(\phi)\frac{\partial\phi^i}{\partial x_\alpha}\frac{\partial\phi^j}{\partial x_\beta}
+\langle\psi, D\hspace{-0.25cm}/\hspace{0.1cm}\psi\rangle\right]\sqrt{g}\,dx.
\eeq

Critical points of $L(\phi,\psi)$ are called Dirac-harmonic maps, which are natural extensions of
harmonic maps and harmonic spinors. In fact, when $\psi=0$, $L(\phi,0)=\int_M |d\phi|^2\,dv_g$ is
the Dirichlet energy functional of $\phi:M\to N$, and its critical points are harmonic maps that
have been extensively studied (see Lin-Wang \cite{lin-wang} for relevant references).  On the other hand, when
$\phi=\hbox{constant}:M\to N$ is a constant map, $L(\hbox{constant},\psi)=\int_M\langle\psi,D\hspace{-0.25cm}/\hspace{0.1cm}\psi\rangle\,dv_g$
is the Dirac functional of $\psi\in (\Gamma\Sigma)^k$, and its critical points are harmonic spinors $\partial\hspace{-0.2cm}/\hspace{0.1cm}\psi=0$
that have also been well studied (see Lawson-Michelsohn \cite{lawson-michelsohn}).

One of the main interests we have is to study regularity of weakly Dirac-harmonic maps. To do it,
we introduce the natural Sobolev space in which the functional $L(\cdot,\cdot)$ is well defined.
Recall the Sobolev space $H^1(M,N)$ is defined by
$H^1(M,N)=\left\{u\in H^1(M,\mathbb R^K): \ u(x)\in N \ \hbox{ a.e. }x\in M\right\}.$
\bde{}\label{sobolev-spinor} For $\phi\in H^1(M,N)$, the set of sections $\psi\in\Gamma(\Sigma\tensor \phi^*TN)$
is defined to be all $\psi=(\psi^1,\cdots,\psi^K)\in \left(\Gamma\Sigma\right)^K$ such that
$$\sum_{i=1}^K \nu_i \psi^i(x)=0 \ \hbox{ a.e. }x\in M ,\ \forall \nu=(\nu_1,\cdots, \nu_K)\in (T_{\phi(x)}N)^\perp.$$
We say that $\psi=(\psi^1,\cdots,\psi^K)\in S^{1,\frac43}\left(\Gamma(\Sigma\tensor\phi^*TN)\right)$
if $d\psi^i\in L^{\frac43}(M)$ and $\psi^i\in L^4(M)$ for all $1\le i\le K$.
\ede
\bde{} \label{weak-dirac-harm}
A pair of maps $\left(\phi,\psi\right)\in H^1(M,N)\times S^{1,\frac43}\left(\Gamma(\Sigma\tensor\phi^*TN)\right)$
is called a weakly Dirac-harmonic map, if it is a critical point of $L\left(\cdot,\cdot\right)$ over
the Sobolev space $H^1(M,N)\times S^{1,\frac43}\left(\Gamma(\Sigma\tensor\phi^*TN)\right)$.
\ede
\bre{} First, the H\"older inequality implies that
if $\left(\phi,\psi\right)\in H^1(M,N)\times S^{1,\frac43}\left(\Gamma(\Sigma\tensor\phi^*TN)\right)$, then
$$\left\vert\int_M\left\langle \psi,D\hspace{-0.25cm}/\hspace{0.1cm}\psi\right\rangle \,dv_g\right\vert
\leq C\left\Vert \psi\right\Vert _{L^4(M)}\left[\left\Vert d
\psi\right\Vert_{L^{\frac43}(M)}+\left\Vert d\phi\right\Vert_{L^2(M)}\left\Vert \psi\right\Vert _{L^4(M)}\right]
<+\infty.$$
Hence $L\left(\phi,\psi\right)$ is well defined for any $(\phi,\psi)\in H^1(M,N)\times S^{1,\frac43}\left(\Gamma(\Sigma\tensor\phi^*TN)\right)$.
Second, it is useful to notice that
\beq{}\int_M\left\langle \psi,D\hspace{-0.25cm}/\hspace{0.1cm}\psi\right\rangle \,dv_g
=\int_M\hbox{Re}\left\langle \psi,D\hspace{-0.25cm}/\hspace{0.1cm}\psi\right\rangle \,dv_g, \label{self_ad}
\eeq
where $\hbox{Re}(z)$ denotes the real part for $z\in \mathbb C$. In fact, since
$D\hspace{-0.25cm}/\hspace{0.1cm}$ is self-adjoint, i.e.
$$\int_M\left\langle \psi,D\hspace{-0.25cm}/\hspace{0.1cm}\xi\right\rangle \,dv_g
=\int_M\left\langle D\hspace{-0.25cm}/\hspace{0.1cm}\psi,\xi\right\rangle \,dv_g, \ \forall \ \psi,
\ \xi\in S^{1,\frac43}\left(\Gamma(\Sigma\tensor\phi^*TN)\right),$$
we have
$$\int_M \overline{\langle\psi,D\hspace{-0.25cm}/\hspace{0.1cm}\psi\rangle}\,dv_g
=\int_M \langle D\hspace{-0.25cm}/\hspace{0.1cm}\psi,\psi\rangle\,dv_g=\int_M \langle\psi,D\hspace{-0.25cm}/\hspace{0.1cm}\psi\rangle\,dv_g.$$
This yields (\ref{self_ad}).
\ere

The Euler-Lagrange equation of a Dirac-harmonic maps $(\phi,\psi)\in H^1(M,N)\times S^{1,\frac43}\left(\Gamma(\Sigma\tensor\phi^*TN)\right)$
is (cf. \cite{chen-jost-li-wang1}):
\begin{eqnarray}
\tau(\phi)&=&\mathcal R^N(\phi,\psi),\label{dirac_harm_eqn}\\
D\hspace{-0.25cm}/\hspace{0.1cm}\psi&=&0, \label{dirac_eqn}
\end{eqnarray}
where $\tau(\phi)$ is the tension field of $\phi$ given by
$$\tau(\phi)=\hbox{tr}(\nabla^{M\tensor\phi^*TN} d\phi)=\left(\Delta\phi^i+g^{\alpha\beta}\Gamma_{jl}^i(\phi)
\frac{\partial\phi^j}{\partial x_\alpha}\frac{\partial\phi^l}{\partial x_\beta}\right)\frac{\partial}{\partial y_i}(\phi),$$
and $\mathcal R^N(\phi,\psi)\in \Gamma(\phi^*TN)$ is defined by
$$\mathcal R^N(\phi,\psi)
=\frac{1}{2}\sum R_{lij}^m(\phi)\left\langle \psi^i,\nabla \phi^l\circ
\psi^j\right\rangle \frac{\partial}{\partial y_m}(\phi).$$
Here $\Gamma_{jl}^i(\phi)$ is the Christoffel symbol of the Levi-Civita connection of $N$,
$\nabla\phi^l\circ \psi^j$  denotes the Clifford multiplication of the vector field
$\nabla\phi^l$ with the spinor $\psi^j$,  and $R^m_{lij}$ is a component of
the Riemannian curvature tensor of $(N,h)$.

Among other things, Chen-Jost-Li-Wang proved in \cite{chen-jost-li-wang2} (Theorem 2.2 and 2.3)
that {\it if $(M^2,g)$ is a spin Riemann surface and $N=S^{K-1}\subset\mathbb R^K$ is the standard sphere,
then any weakly Dirac-harmonic map $(\phi,\psi)$ is in
$C^\infty(M^2,S^{K-1})\times C^\infty\left(\Gamma(\Sigma\times\phi^*(TS^{K-1}))\right)$}, which was
extended to any compact hypersurface in $\mathbb R^K$ by Zhu \cite{zhu}.
The crucial observation in \cite{chen-jost-li-wang2} is that the nonlinearity 
in (\ref{dirac_harm_eqn}) is of Jacobian determinant structure.
This theorem is an extension of that on harmonic maps by H\'elein \cite{helien}.
Our motivation in this paper is: (1) extend the above theorem to all Riemannian manifold $N\subset\mathbb R^K$,
and (2) study the regularity problem of stationary Dirac-harmonic maps in higher dimensions $n\ge 3$.

Denote by ${\it i}_M>0$ the injectivity radius of $M$. For $0<r<{\it i}_M$ and $x\in M$, denote by $B_r(x)$
the geodesic ball in $M$ with center $x$ and radius $r$.
Our first result is an $\epsilon$-regularity theorem.
\bth{}\label{small_reg} For $n\ge 2$,
there exists $\varepsilon _{0}>0$ depending only on $(M,g)$ and $(N,h)$ such that if $(\phi,\psi)\in H^1(M,N)\times S^{1,\frac43}\left(\Gamma(\Sigma\tensor\phi^*TN)\right)$ is a weakly Dirac-harmonic map satisfying, for some $x_0\in M$
and $0<r_0\le \frac12{\it i}_M$,
\beq{}\label{small_reg1}
\underset{x\in B_{r_0}(x_0), 0<r\le r_0}{\sup }\left\{ \frac{1}{r^{n-2}}\int_{B_{r}\left(
x\right) }\left( \left\vert d\phi\right\vert ^{2}+\left\vert \psi
\right\vert ^{4}\right) dv_g\right\} <\varepsilon_0^2,
\eeq
then $(\phi,\psi)$ is smooth in $B_{r_0}(x_0)$.
\eth

Since $\int_M (|d\phi|^2+|\psi|^4)$ is conformally invariant when $n(=\hbox{dim}\ M)=2$ (see \cite{chen-jost-li-wang1} Lemma 3.1),
it is not hard to see that there exists $0<r_0=r_0(M, \phi,\psi)\le {\it i}_M$ such that
$$\underset{x\in M}\sup \int_{B_{r_0}(x)}\left(|d\phi|^2+|\psi|^4\right)\le \epsilon_0^2,$$
where $\epsilon_0>0$ is the same constant as in theorem \ref{small_reg}. Hence, as an immediate consequence of
theorem \ref{small_reg}, we have
\bth{}\label{2d_reg} For $n=2$, assume $\left(\phi,\psi\right)\in H^1(M,N)\times S^{1,\frac43}\left(\Gamma(\Sigma\tensor\phi^*TN)\right)$
is a weakly Dirac-harmonic map. Then $\left(\phi,\psi\right)\in C^\infty(M, N)\times C^\infty(\Gamma(\Sigma\tensor\phi^*TN))$.
\eth

We would remark that theorem \ref{2d_reg} has been proved by Chen-Jost-Li-Wang \cite{chen-jost-li-wang2}
when the target manifold $N=S^{K-1}\subset\mathbb R^K$ is the standard sphere.

For $n\ge 3$, it is well-known in the context of harmonic maps that in order for a harmonic map to enjoy partial regularity,
we need to pose the stationarity condition (see, e.g. Evans \cite{evans}, Bethuel \cite{bethuel}, and Rivier\'e \cite{riviere}).

For the same purpose, we also introduce the notion of stationary Dirac-harmonic maps.
\bde{}\label{stat-dirac} We call a  weakly Dirac-harmonic map $$(\phi,\psi)\in H^1\left(M,N\right) \times \mathcal S^{1,\frac43}
\left(\Gamma\left(\Sigma \otimes \phi^*TN\right)\right)$$
to be  a stationary Dirac-harmonic map, if, in additions, it is a critical point of
$L\left(\phi,\psi\right)$ with respect to the domain variations, i.e., for any family of
diffeomorphisms $F_t(x):=F(t,x)\in C^1((-1,1)\times M, M)$
with $F_0(x) =x$ for $x\in M$, and $F_t(x) =x$
for any $x\in \partial M$ and $t\in(-1,1) $ when $\partial M\neq\emptyset$, then we have
\beq{}\label{stat_dirac1}
\left.\frac{d}{dt}\right\vert _{t=0}\left[ \int_M\left(|d\phi_{t}|^2+\left\langle \psi _{t},D\hspace{-0.25cm}/\hspace{0.1cm}\psi _{t}\right\rangle\right) dv_{g}
\right] =0,
\eeq
where $\phi_t(x)=\phi(F_t(x))$ and $\psi_{t}=\psi( F_{t}(x)) $.
\ede

Motivated by \cite{chen-jost-li-wang1}, we define {\it stress-energy tensor} $\mathcal S$
for a stationary Dirac-harmonic map $(\phi,\psi)$ by
\beq{}\label{stress_tensor}
\mathcal S_{\alpha\beta}=\left\langle\frac{\partial\phi}{\partial x_\alpha}, \frac{\partial\phi}{\partial x_\beta}\right\rangle
-\frac12|d\phi|^2\delta_{\alpha\beta}
+\frac12\hbox{Re}\left\langle \psi, \frac{\partial}{\partial x_\alpha}\circ\nabla_{\frac{\partial}{\partial x_\beta}}\psi\right\rangle,
\ 1\le\alpha,\beta\le n,
\eeq
where $\{\frac{\partial}{\partial x_\alpha}\}$ is a local coordinate frame on $M$.

It turns out that the stationarity property is equivalent to that the stress-energy tensor $\mathcal S$ is divergence free
(see Lemma \ref{free-stress}):
\beq{}\label{stationary_id1}
\sum_{\alpha, \beta}\frac{\partial}{\partial x_\alpha}\left(\sqrt {g} g^{\alpha\beta}\mathcal S_{\beta\gamma}\right)=0,
\ 1\le\gamma\le n, \eeq
in the sense of distributions.

An immediate consequence of (\ref{stationary_id1}), we prove in \S4 the following Liouville property of stationary Dirac-harmonic maps.
\bth{} \label{liouville} For $n\ge 3$, let $(M,g)=(\mathbb R^n, g_0)$ be the $n$-dimensional Euclidean space associated with the spinor bundle
$\Sigma$. If $(\phi,\psi)\in H^1(\mathbb R^n,N)\times \mathcal S^{1,\frac43}(\Gamma(\Sigma\tensor \phi^*TN))$
is a stationary Dirac-harmonic map, then $\phi\equiv\hbox{constant}$ and $\psi\equiv 0$.
\eth

In dimensions $n\ge 3$, the stationarity property is a necessary condition for smoothness of weakly Dirac-harmonic maps.
In fact, Chen-Jost-Li-Wang \cite{chen-jost-li-wang1} proved that any smooth Dirac-harmonic map $(\phi,\psi)\in C^\infty(M,N)\times C^\infty(\Gamma(\Sigma\times\phi^*TN))$ has its stress-energy tensor divergence free and hence is a stationary Dirac-harmonic map. 
Hence theorem \ref{liouville} extends a corresponding Liouville theorem on smooth Dirac-harmonic maps by Chen-Jost-Wang \cite{chen-jost-wang}.

An important implication of (\ref{stationary_id1}) is the following monotonicity inequality (see \S4 below):
{\it there exist $0<r_0<{\it i}_M$ and $C_0>0$ depending only on $(M,g)$ such that if
$(\phi,\psi)\in H^1\left(M,N\right) \times \mathcal S^{1,\frac43}\left(\Gamma\left(\Sigma \otimes \phi^*TN\right)\right)$ is
a stationary Dirac-harmonic map, then for any $x_0\in M$ and $0<r\le r_0$, it holds
\begin{eqnarray}\label{mono_id012}
\frac{d}{dr}\left(e^{C_0r}r^{2-n}\int_{B_r(x)}\left\vert d\phi\right\vert
^{2}dv_g\right) &\ge& e^{C_0r}r^{2-n}\int_{\partial B_r(x)}
2\vert \frac{\partial \phi}{\partial r}\vert ^{2}\,dH^{n-1}\\
&&+e^{C_0r}r^{2-n}\int_{\partial B_r(x_0)}{\rm{Re}}\langle \psi ,\frac{\partial }
{\partial r}\circ \nabla _{\frac{\partial }{\partial r}}\psi\rangle\,dH^{n-1}. \nonumber
\end{eqnarray}
}

However, we should point out that (\ref{mono_id012}) doesn't yield the renormalized energy
$e^{C_0r}r^{2-n}\int_{B_r(x)}\left\vert d\phi\right\vert
^{2}dv_g$
is monotone increasing with respect to $r$, since the second term of the right hand side of (\ref{mono_id012})
$$e^{C_0r}r^{2-n}\int_{\partial B_r(x)}\hbox{Re}\langle \psi ,
\frac{\partial }{\partial r}\circ \nabla _{\frac{\partial }{\partial r}}\psi\rangle
dH^{n-1}$$
may change signs. In order to utilize (\ref{mono_id012}) to control
$r^{2-n}\int_{B_r(x)}\left(|d\phi|^2+|\psi|^4\right)\,dv_g,$
we need to assume $d\psi\in L^p$ for some $p>\frac{2n}3$.  In fact, we have
\bth{}
\label{nd_reg} For $n\ge 3$, let $(\phi,\psi)\in H^1\left(M,N\right) \times \mathcal S^{1,\frac43}
\left(\Gamma\left(\Sigma \otimes \phi^*TN\right)\right)$ be a stationary Dirac-harmonic map. If, in additions,
$d\psi\in L^p(M)$ for some $p>\frac{2n}3$, then there exists a closed subset $\mathcal S(\phi)\subset M$,
with $H^{n-2}(\mathcal S(\phi))=0$, such that $(\phi,\psi)\in C^\infty(M\setminus \mathcal S(\phi))$.
\eth

Now let's outline the main ingredients to prove theorem \ref{small_reg} as follows. \\
(1) We observe that the Dirac-harmonic property is
invariant under totally geodesic, isometric embedding. More precisely, let $\Phi: (N,h)\to(\tilde N,\tilde h)$
be a totally geodesic, isometric embedding map. If, for $\phi:M\to N$ and $\psi\in \Gamma(\Sigma\tensor\phi^*TN)$,
$(\phi,\psi)$ is a weakly Dirac-harmonic map, then for $\tilde \phi=\Phi(\phi):M\to\tilde N$ and
$\tilde\psi=\Phi_*(\psi)=\psi^i\tensor \frac{\partial}{\partial z_i}(\tilde\phi)\in\Gamma(\Sigma\tensor {\tilde\phi}^*T\tilde N)$, $(\tilde\phi,\tilde\psi)$ is a weakly Dirac-harmonic map. \\
(2) By employing the enlargement argument by H\'elein \cite{helein1, helein2} in the context of harmonic maps, we can assume that $TN\big|_{\phi(M)}$
is trivial so that there exists an orthonormal tangent frame $\{e_i\}_{i=1}^k$ on $\phi^*TN$.\\
(3) We use this moving frame to rewrite the Dirac-harmonic map equation (\ref{dirac_harm_eqn}) into the form
\beq{} \label{new_dirac_harm_eqn} d^*\left(\langle d\phi,e_i\rangle\right)=\sum_{j}\Theta_{ij}\langle d\phi, e_j\rangle,
\eeq
where $\Theta=(\Theta_{ij})\in L^2\left(B_{r_0}(x_0), {\it so}(n)\tensor \wedge^1(\mathbb R^n)\right)$ satisfies
$|\Theta|\le C(|d\phi|+|\psi|^2)$. \\
(4) The smallness condition (\ref{small_reg1}) guarantees that
we can apply the Coulomb gauge construction, due to Revier\'e \cite{riviere1} ($n=2$) and Rivier\'e-Struwe \cite{riviere-struwe} ($n\ge 3$),
to further rewrite (\ref{new_dirac_harm_eqn}) into an equation in which the nonlinearity has the Jacobian determinant structure similar to
that of harmonic maps. \\
(5) We utilize the duality between the Hardy space and BMO space to obtain an decay estimate in the Morrey space, which
yields the H\"older continuity of $\phi$. \\
(6) By adapting the hole-filling technique developed by Giaquinta-Hildebrandt \cite{giaquinta-hildebrandt}
in the context of harmonic maps, we establish the higher order regularity of $(\phi,\psi)$. We point out that in dimension two,
a different proof of higher order regularity of Dirac-harmonic maps has been provided by Chen-Jost-Li-Wang \cite{chen-jost-li-wang2}.

As a byproduct of the rewriting of Dirac-harmonic maps under the above Coulomb gauge frame, we also obtain a
convergence theorem of weakly convergent sequences of approximate Dirac-harmonic maps in dimension two,
which extends a corresponding convergence of approximate harmonic maps from surfaces due to Bethuel \cite{bethuel1}
(see also Freire-M\"uller-Struwe \cite{FMS}, Wang \cite{wang1}, and Rivier\'e \cite{riviere1}).
More precisely, we have
\bth{}\label{conv_dirac_hm}
For $n=2$, let $\left(\phi_p,\psi_p\right)\in H^1\left(M,N\right) \times \mathcal S^{1,\frac43}
\left(\Gamma\left(\Sigma \otimes \phi^*TN\right)\right)$ be a sequence of weak solutions to
the approximate Dirac-harmonic map equation
\begin{eqnarray}
\tau\left(\phi_p\right)&=&\mathcal R^N(\phi_p,\psi_p)+u_p \label{approx_dirac1}\\
D\hspace{-0.25cm}/\hspace{0.1cm}\psi_p&=&v_p. \label{approx_dirac1}
\end{eqnarray}
Assume that $u_p\rightarrow 0$ strongly in $H^{-1}(M)$ and $v_p\rightharpoonup 0$ weakly in $L^{\frac43}(M)$.
If $\phi_p\rightharpoonup\phi$ in $H^1(M,N)$ and $\psi_p\rightharpoonup\psi$ in $\mathcal S^{1,\frac43}$,
then $(\phi,\psi)\in H^1\left(M,N\right) \times \mathcal S^{1,\frac43}
\left(\Gamma\left(\Sigma \otimes \phi^*TN\right)\right)$ is a weakly Dirac-harmonic map.
\eth
The paper is organized as follows. In \S2, we rewrite the equation of Dirac-harmonic maps via moving frames.
In \S3, we use the Coulomb gauge construction, duality between Hardy space and BMO space,
and a decay estimate in Morrey space to first prove the H\"older continuity part of Theorem \ref{small_reg},
and then adopt the hole-filling technique by Giaquinta-Hildebrandt \cite{giaquinta-hildebrandt}
to prove the higher order regularity part of Theorem \ref{small_reg}.
In \S4, we discuss various properties of stationary Diac-harmonic maps and prove
Theorem \ref{liouville} and Theorem \ref{nd_reg}. In \S5, we prove the convergence Theorem \ref{conv_dirac_hm}.

\section{Dirac-harmonic maps via moving frames}

In this section, we first show that a Dirac harmonic map $(\phi,\psi)$ is invariant under
a totally geodesic, isometric embedding so that H\'elein's enlargement argument (cf. \cite{helein1, helein2})
guarantees that we can assume there is an orthonormal frame $\{e_i\}_{i=1}^k$ of $\phi^*TN$.
Then employing this orthonormal frame we write the equation of Dirac-harmonic maps
into the form (\ref{new_dirac_harm_eqn}).

We begin with
\bpr{}\label{invariance} Let $(\widetilde N,\widetilde h)$ be another compact Riemannian manifold without boundary
and $f:\left(N,h\right)\to \left(\widetilde N,\widetilde h\right)$ be a totally geodesic,
isometric embedding. Let $(\phi, \psi)\in H^1(M,N)\times \mathcal S^{1,\frac43}
\left(\Gamma(\Sigma M\otimes \phi^*TN)\right)$ be a weakly Dirac-harmonic map.
Define $\widetilde u=f(u)\in H^1(M, \widetilde N)$ and
$\widetilde\psi=f_*(\psi)=\psi^i\tensor \frac{\partial}{\partial z_i}(\widetilde\phi)
\in \mathcal S^{1,\frac43}(\Gamma(\Sigma\tensor(\widetilde\phi)^*T\widetilde N))$.
Then $\left(\widetilde\phi,\widetilde{\psi}\right)$ is also a weakly Dirac-harmonic map.
\epr
\proof By the chain rule formula of tension fields (cf. Jost \cite{jost}),  we have
$$
\tau\left(\widetilde\phi\right)=\hbox{tr}\left[\nabla^{f*T\widetilde N}df
\left(d\phi, d\phi\right)\right]+f_{\ast}\left(\tau\left(\phi\right)\right)
=f_{\ast}\left(\tau\left(\phi\right)\right)=f_{\ast}\left(\mathcal R^N(\phi,\psi)\right),
$$
where we have used the fact that $f$ is totally geodesic, i.e. $\nabla^{f^*T\widetilde N}df=0$, and
the Dirac-harmonic map equation (\ref{dirac_harm_eqn}).

Set $\widehat N=f(N)$. Then $(\hat N, \widetilde h)$ is a totally geodesic, submanifold of
$(\widetilde N, \widetilde h)$. Moreover,  if $y=(y_1,\cdots, y_k)$ is a local coordinate system on $N$,
then $z=(z_1,\cdots,z_k)=f(y)$ is a local coordinate system on $\widehat N$ and
$\frac{\partial}{\partial z_i}=f_*(\frac{\partial}{\partial y_i}), 1\le i\le k$ is a local coordinate frame
on $\widehat N$. Since $f:(N,h)\to (\widehat N,\widetilde h)$ is an isometry, we have
\begin{eqnarray*}f_*(\mathcal R^N(\phi,\psi))
&=&f_*\left(\frac12 (R^N)^m_{lij}(\phi)\langle\psi^i,\nabla\phi^l\circ\psi^j\rangle \frac{\partial}{\partial y_m}(\phi)\right)\\
&=&\frac12(R^{\widehat N})^m_{lij}(\widetilde\phi)\langle\psi^i,\nabla\widetilde\phi^l\circ\psi^j\rangle
\frac{\partial}{\partial z_m}(\widetilde\phi)\\
&=&\mathcal R^{\widehat N}(\widetilde\phi,\widetilde\psi)=\mathcal R^{\widetilde N}(\widetilde\phi,\widetilde\psi),
\end{eqnarray*}
where we have used the fact that $(R^{\widehat N})^m_{lij}(\widetilde \phi)=(R^{\widetilde N})^m_{lij}(\widetilde\phi)$
in the last two steps, which follows from the Gauss-Codazzi equation since $\widehat N\subseteq\widetilde N$ is
a totally geodesic submanifold.

To see that $\widetilde{\psi}$ satisfies the equation (\ref{dirac_eqn}), denote $\widetilde{D\hspace{-0.25cm}/\hspace{0.1cm}}$
as the Dirac operator along the map $\widetilde\phi$. Then it follows from \cite{chen-jost-li-wang1} (2.6)
that
$$
\widetilde{D\hspace{-0.25cm}/\hspace{0.1cm}}\widetilde{\psi}=f_{\ast}\left(  D\hspace{-0.25cm}/\hspace{0.1cm}\psi\right)
+(\nabla\phi^i\circ\psi^j)\otimes\nabla^{f^*T\widetilde N}df\left(\frac{\partial}{\partial y_i},
 \frac{\partial}{\partial y_j}\right)=0,
$$
where we have used the fact that both $D\hspace{-0.25cm}/\hspace{0.1cm}\psi=0$ and $\nabla^{f^*T\widetilde N}df=0$.
\qed\\

With the help of Proposition \ref{invariance}, we can now adapt the same enlargement argument
as that by H\'elein \cite{helein1, helein2} and assume that $(N,h)$ is parallelized.
Hence there exist a global orthonormal frame $\left\{\hat e_i\right\}_{i=1}^{k}$ on $(N,h)$.
Set $e_i(x)=\hat e_i(\phi(x))$, $1\le i\le k$. Then $\{e_i\}$ is an orthonormal frame along
$\phi^*TN$. Using this frame, we can write the spinor field $\psi$ along map $\phi$ as
$$
\psi=\sum_{i=1}^{k}\psi^{i}\otimes e_{i},\ \psi^i\in\Gamma(\Sigma), \ 1\le i\le k.$$

Let $\{\frac{\partial}{\partial x_\alpha}\}_{\alpha=1}^n$ be a local coordinate frame on $M$.
Recall the tension field of $\phi$ is defined by (cf. Jost \cite {jost}):
$$
\tau\left(\phi\right)  =g^{\alpha\beta}\nabla_{\frac{\partial}{\partial x_\alpha}}^{\phi^*TN}
\left(\frac{\partial\phi}{\partial x_\beta}\right).$$
Denote the components of $\tau(\phi)$ and $\mathcal D\psi$ with respect to the frame $\{e_i\}$ by
$$\tau^{i}\left(\phi\right)=\left\langle \tau\left(\phi\right),e_{i}\right\rangle_{\phi^*h}, \ 1\leq i\leq k, $$
$$(D\hspace{-0.25cm}/\hspace{0.1cm}\psi)^{i}=\left\langle D\hspace{-0.25cm}/\hspace{0.1cm}\psi,e_i\right\rangle_{\phi^*h}, \ 1\leq i\leq k.$$

Under these notations, we have
\ble{}\label{dirac_hm_eqn} If $(\phi,\psi)\in H^1\left(M,N\right) \times \mathcal S^{1,\frac43}
\left(\Gamma\left(\Sigma \otimes \phi^*TN\right)\right)$ is a weakly Dirac-harmonic map,  then it holds, for
$1\leq i\leq n$,
\begin{eqnarray}
\left(D\hspace{-0.25cm}/\hspace{0.1cm}\psi\right)^i&  =&0\label{dirac_hm1}\\
\tau^i\left(\phi\right)&=& R^{N}(\phi)\left( e_i,e_j,e_l,e_m\right)
\left\langle \phi_*(\xi_\alpha), \ e_j\right\rangle _{\phi^*h}\left
\langle \psi^m,\ \xi_\alpha\circ\psi^l\right\rangle,
\label{dirac-hm2}
\end{eqnarray}
where $\{\xi_\alpha\}_{\alpha=1}^n$ is a local orthonormal frame on $M$.
\ele
\proof It suffices to prove (\ref{dirac-hm2}).  To do this,
let $\left\{\phi_t\right\}$ be a variation of $\phi$ such
that $\frac{\partial\phi_t}{\partial t}\big|_{t=0}=\eta=\eta^{i}e_i$ for
$(\eta^1,\cdots, \eta^k)\in C^\infty_0(M,\mathbb R^k)$.
Then we have $\psi_t=\psi^i\tensor e_i(\phi_t)$.
Then we have
\begin{eqnarray*}
& & \frac{\partial}{\partial t}D\hspace{-0.25cm}/\hspace{0.1cm}\psi_t\\
&=&\partial\hspace{-0.2cm}/\hspace{0.1cm}\psi^i\otimes\nabla_{\frac{\partial}{\partial t}}e_i(\phi_t)
+(\xi_\alpha\circ \psi^i) \otimes
\nabla_{\frac{\partial}{\partial t}}\nabla_{\xi_\alpha}e_i(\phi_t)\\
&=&\partial\hspace{-0.2cm}/\hspace{0.1cm}\psi^{i}\otimes\nabla_{\frac{\partial}{\partial t}}e_i(\phi_t)
+(\x_\alpha\circ\psi^i)\otimes
\nabla_{\xi_\alpha}\nabla_{\frac{\partial}{\partial t}}e_i(\phi_t)\\
&&+(\xi_\alpha\circ\psi^{i})\otimes R^{\phi_t^*TN}(\phi_t)
(\frac{\partial}{\partial t},\xi_\alpha)e_{i}(\phi_t)\\
&=&D\hspace{-0.25cm}/\hspace{0.1cm}(\psi^i\otimes\nabla_{\frac{\partial}{\partial t}}e_i(\phi_t))
+(\xi_\alpha\circ\psi^i)\otimes R^{\phi_t^*TN}(\phi_t)(\frac{\partial}{\partial t},
\xi_\alpha)  e_i(\phi_t).
\end{eqnarray*}
This, combined with the fact that $D\hspace{-0.25cm}/\hspace{0.1cm}\psi=0$ and $D\hspace{-0.25cm}/\hspace{0.1cm}$ is self-adjoint, implies
\begin{eqnarray*}
& & \frac{d}{dt}|_{t=0}\int_{M}\left\langle \psi_t,D\hspace{-0.25cm}/\hspace{0.1cm}\psi_t\right\rangle \,dv_g\\
&=&\int_M\left \langle \frac{\partial}{\partial t}\big|_{t=0}\psi_t, D\hspace{-0.25cm}/\hspace{0.1cm}\psi\right\rangle
+\int_M \left\langle\psi, \frac{\partial}{\partial t}\big|_{t=0}D\hspace{-0.25cm}/\hspace{0.1cm}\psi_t\right\rangle\\
&=& \int_{M}\left\langle \psi^i,\ \xi_\alpha\circ \psi^j\right\rangle
\left\langle e_i,\ R^{\phi^*TN}(\phi)(\eta,\phi_*(\xi_\alpha))e_j\right\rangle \,dv_g\\
&=&\int_{M}\eta^l\langle\phi_*(\xi_\alpha),\ e_m\rangle
\left\langle \psi^i,\ \xi_\alpha\circ\psi^j\right\rangle
\left\langle e_i,\ R^{N}\left(e_l,e_m\right)e_j\right\rangle\,dv_g,
\end{eqnarray*}
On the other hand, it is well-known that
$$\frac{d}{dt}\big|_{t=0}\int_M |d\phi_t|^2\,dv_g
=2\int_M \langle \tau(\phi), e_l\rangle \eta^l\,dv_g=2\int_M\tau^l(\phi)\eta^l\,dv_g.$$
Hence, combining these formula together, we obtain (\ref{dirac-hm2}). \qed

\section{The $\epsilon$-decay estimate and regularity theorem}

In this section, we utilize the skew-symmetry of the nonlinearity in the right hand side
of the Dirac-harmonic map equation (\ref{dirac-hm2}) and adapt the Coulomb gauge construction
technique developed by Rivier\'e \cite{riviere1} ($n=2$) and Rivier\'e-Struwe \cite{riviere-struwe}
($n\ge 3$) to establish an energy decay estimate for Dirac-harmonic maps in Morrey spaces under
the smallness condition.  As consequences, we prove the H\"older continuity part of
Theorem \ref{small_reg} and Theorem \ref{2d_reg}.

Since the regularity issue is a local result, we assume, for simplicity of presentation, that
for $x_0\in M$, the geodesic ball $B_{{\it i}_M}(x_0)\subset M$ with the metric $g$
is identified by $(B_2,g_0)$. Here $B_2$ is the ball centered at $0$ and
radius $2$ in $\mathbb R^n$, and $g_0$ is the Euclidean metric on $\mathbb R^n$.
We also assume that the spin bundle $\Sigma$ restricted in $B_2$ is given by
$\Sigma\big|_{B_2}\equiv B_2\times \mathbb C^L$, with $L={\rm{rank}}_{\mathbb C}\Sigma$.

Let $(\phi,\psi)\in H^1(B_2, N)\times \mathcal S^{1,\frac43}(B_2, \mathbb C^L\tensor\phi^*TN)$ be a weakly Dirac-harmonic map,
and $\{e_i\}_{i=1}^k$ be an orthonormal frame of $\phi^*TN$ given as in \S2. Write
$\psi=\psi^i\tensor e_i$ for some $\psi^i\in\mathbb C^L$, $1\le i\le k$.

Now we define $\Omega$, the $k\times k$ matrix whose entries are $1$-forms, by
\beq{}\label{1-form}
\Omega_{ij}=\sum_{\alpha=1}^n\left[\sum_{l,m=1}^k R^{N}(\phi)
\left(e_i, e_j, e_l, e_m\right)\langle \psi^m,\frac{\partial}{\partial x_\alpha}\circ\psi^l\rangle\right]
\,dx_\alpha, \text{, \ for }1\leq i, j\leq k.
\eeq
Then we have the following simple fact.
\bpr{}\label{skew-sym} Let $\Omega$ be given by (\ref{1-form}). Then
$\Omega_{ij}$ is real valued for any $1\leq i, j\leq k$, and
$\Omega$ is skew-symmetric, i.e.
$$\Omega_{ij}=-\Omega_{ji}, \ 1\le i, j\le k.
$$
\epr
\proof First observe that the skew-symmetry of Clifford multiplication $\circ$
and the properties of Hermitian metric $\left\langle \cdot,\cdot\right\rangle$
give
$$\overline{\left\langle \psi^m, \frac{\partial}{\partial x_\alpha}\circ
\psi^l\right\rangle}=\left\langle \frac{\partial}{\partial x_\alpha}\circ
\psi^l, \psi^m\right\rangle
=-\left\langle \psi^l, \frac{\partial}{\partial x_\alpha}\circ\psi^m\right\rangle.$$
On the other hand,  the curvature operator $R^{N}(\phi)(\cdot,\cdot,\cdot,\cdot)$ is skew-symmetric
in its last two components:
$$ R^{N}(\phi)\left(\cdot, \cdot, e_l, e_m\right)= -R^{N}(\phi)\left(\cdot, \cdot, e_m, e_l\right).$$
Thus we conclude that
$$\overline{\Omega_{ij}}=\Omega_{ij}$$
so that $\Omega$ is real valued.
$\Omega_{ij}=-\Omega_{ji}$ follows from skew-symmetry of $R^{N}(\phi)(\cdot,\cdot,\cdot,\cdot)$
with respect to its first two components. \qed\\

In terms of $\Omega$,  the Dirac-harmonic map equation (\ref{dirac-hm2}) can be written as
\beq{}\label{dirac-hm3}
\tau^i\left(\phi\right)  =\sum_{j=1}^k\Omega_{ij}\cdot\left\langle d\phi, e_j\right\rangle,
\text{ \ }1\leq i\leq k,
\eeq
where $\cdot$ denotes the inner product of 1-forms, and
$
\left\langle d\phi, e_j\right\rangle =\sum_{\alpha=1}^{n}\left\langle
\frac{\partial \phi}{\partial x_\alpha}, e_{j}\right\rangle\,dx_\alpha.$

Denote by $d^*$ the conjugate operator of $d$. Then we have, for $1\le i\le k$,
$$d^*\left(\langle d\phi,e_i\rangle\right)
=\langle\tau(\phi),e_i\rangle+\langle d\phi, de_i\rangle
=\tau^i(\phi)+\langle de_i, e_j\rangle\cdot \langle d\phi, e_j\rangle.$$
Hence we have
\beq{}\label{dirac-hm4}
d^*\left(\langle d\phi,e_i\rangle\right)=\sum_{l=1}^k \Theta_{il}\cdot \langle d\phi,e_l\rangle;
 \ \Theta_{ij}\equiv\Omega_{ij}+\langle de_i, e_j\rangle, \ \forall 1\le i, j\le k.
\eeq

Before proving Theorem \ref{small_reg}, we recall the definition of Morrey spaces.
\bde{}\label{morrey} For $1\le p\le n$, $0<\lambda\le n$, and a domain $U\subseteq\mathbb R^n$, the
Morrey space $M^{p,\lambda}(U)$ is defined by
$$
M^{p,\lambda}(U)
:=\left\{f\in L^p_{\hbox{loc}}(U): \|f\|_{M^{p,\lambda}(U)}<+\infty\right\},$$
where
$$\left\Vert f\right\Vert _{M^{p,\lambda}(U)}^p=\sup\left\{r^{\lambda-n}\int_{B_r}|f|^p:\ B_r\subseteq U\right\}.$$
It is easy to see that for $1\leq p\leq n$, $M^{p,n}(U)=L^p(U)$ and  $M^{p,p}(U)$ behaves like
$L^{n}(U)$ from the view of scalings.
\ede

Now we recall the Coulomb gauge construction theorem in Morrey spaces with small Morrey norms,
due to Rivier\'e \cite{riviere1} for $n=2$ and Rivier\'e-Struwe \cite{riviere-struwe} for $n\ge 3$,
which plays a critical role in our proof here.

\ble{}\label{coulomb_gauge}
There exist $\epsilon(n)>0$ and $C(n)>0$
such that if  $R\in L^2(B_1, \rm{so}(k)\tensor \wedge^1 \mathbb R^n)$ satisfies%
\beq{}\label{small_morrey}
\|R\|_{M^{2,2}(B_1)}\le \epsilon(n),
\eeq
then there exist $P\in H^1(B_1, \rm{SO}(k))$ and $\xi\in H^{1}\left(B_1, \rm{so}(k)
\tensor \wedge^2\mathbb R^n\right)$ such that
\begin{eqnarray} \label{coulomb_gauge1}
P^{-1}R P+P^{-1}dP&=&d^{*}\xi \ \hbox{in } B_1 \\
d\xi=0 \ \hbox{ in } B_1, && \xi=0 \ \hbox{ on }\partial B_1.
\end{eqnarray}
Moreover, $\nabla P$ and $\nabla \xi$ belong to $M^{2,2}(B_1)$ with
\beq{}\label{morrey_est}
\|\nabla P\|_{M^{2,2}(B_1)}+\|\nabla\xi\|_{M^{2,2}(B_1)}\le C(n)\|R\|_{M^{2,2}(B_1)}\le C(n)\epsilon(n).
\eeq
Here $\rm{so}\left(k\right)$ denotes Lie algebra of $\rm{SO}\left(k\right)$.
\ele

The crucial step to prove Theorem \ref{small_reg} is the following lemma.
\ble{}\label{morrey_decay} There exist $\epsilon_0>0$ and $\theta_0\in (0,\frac12)$
such that if $(\phi,\psi)\in H^1(B_2,N)\times \mathcal S^{1,\frac43}(B_2, \mathbb C^L\tensor \phi^*TN)$
is a weakly Dirac-harmonic map satisfying
\beq{}\label{small_morrey_norm}
\|\nabla\phi\|_{M^{2,2}(B_2)}^2+\|\psi\|_{M^{4,2}(B_2)}^4\le\epsilon_0^2,
\eeq
then for any $\alpha\in (0,1)$, $\phi\in C^\alpha(B_1,N)$. Moreover,
\beq{}\label{morrey_decay1}
\left[\phi\right]_{C^\alpha(B_1)}\le C\|\nabla\phi\|_{M^{2,2}(B_2)}.
\eeq
\ele
\proof By Proposition \ref{skew-sym} and (\ref{dirac-hm4}), we have
$\Theta=\left(\Theta_{ij}\right)=\left(\Omega_{ij}-\left\langle e_i, de_j\right\rangle \right)
\in L^2\left(B_1, \rm{so}(k)\tensor\wedge^1\mathbb R^n\right)$.
Moreover, (\ref{small_morrey}) implies
$$\|\Theta\|_{M^{2,2}(B_1)}\le C(N)\left[\||\psi|^2\|_{M^{2,2}(B_1)}+\|\nabla\phi\|_{M^{2,2}(B_1)}\right]
\le C(N)\epsilon_0\le\epsilon(n),$$
provide $\epsilon_0>0$ is chosen to be sufficiently small,
where $\epsilon(n)>0$ is the same constant as in Lemma \ref{coulomb_gauge}.
Hence, applying Lemma \ref{coulomb_gauge} with $R$ replaced by $-\Theta$, we conclude that
there are $P\in H^1(B_1, \rm{SO}(k))$ and $\xi\in H^1(B_1, \rm{so}(k)\tensor \wedge^1\mathbb R^n)$
such that
\beq{} \label{coulomb_gauge1}
P^{-1}dP-P^{-1}\Theta P=d^*\xi, \ d\xi =0 \ \hbox{ in }B_1;  \ \xi=0 \ \hbox{ on }\partial B_1,\eeq
and
\beq{}\label{morrey_est1}
\|\nabla P\|_{M^{2,2}(B_1)}+\|\nabla\xi\|_{M^{2,2}(B_1)}\le C(n)\|\Theta\|_{M^{2,2}(B_1)}\le C(n)\epsilon_0.
\eeq
Write $P=(P_{ij})$, $P^{-1}=(P_{ji})$, and $\xi=(\xi_{ij})$. Since $P^{-1}P=I_k$, we have
$dP^{-1}=-P^{-1}dP P^{-1}$.
Multiplying $P^{-1}$ to the equation (\ref{dirac-hm4}) and applying (\ref{coulomb_gauge1}),
we obtain
\begin{eqnarray}
d^*\left[P^{-1}\left(\begin{matrix}\langle d\phi,e_1\rangle \\ \vdots \\
\langle d\phi, e_k\rangle\end{matrix}\right)\right]
&=&\left[dP^{-1}P+P^{-1}\Theta P\right]\cdot
P^{-1}\left(\begin{matrix}\langle d\phi,e_1\rangle \\ \vdots \\ \langle d\phi, e_k\rangle\end{matrix}\right)\nonumber\\
&=&-d^*\xi \cdot
P^{-1}\left(\begin{matrix}\langle d\phi,e_1\rangle\\ \vdots \\ \langle d\phi, e_k\rangle\end{matrix}\right).
\label{dirac-hm5} \end{eqnarray}
The components of (\ref{dirac-hm5}) can be written as
\beq{}\label{dirac-hm6}
-d^*\left(P_{ji}\langle d\phi, e_j\rangle\right)
=d^*\xi_{il}\cdot\left(P_{ml}\langle d\phi, e_m\rangle\right),  \ 1\le i\le k, \ \hbox{ in } B_1.
\eeq

To proceed, recall the definition of BMO spaces. For any domain $U\subseteq\mathbb R^n$, ${\rm{BMO}}(U)$ is
defined to be the set of functions $f\in L^1_{\hbox{loc}}(U)$ such that
$$\left[f\right]_{\hbox{BMO}(U)}\equiv
\sup\left\{\frac{1}{|B_r|}\int_{B_r}|f-\bar f_r|\,dx: \ B_r\subseteq U\right\}<+\infty,$$
where $\bar f_r =\frac{1}{|B_r|}\int_{B_r} f$ is the average of $f$ over $B_r$. By Poincar\'e
inequality, it follows that
\beq{}\label{morrey-bmo}
\left[f\right]_{\hbox{BMO}(U)}\le C\|\nabla f\|_{M^{p,p}(U)}, \ \forall 1\le p\le n.
\eeq

For any $0<R\le \frac12$, let $B_R\subset B_1$ be an arbitrary ball of radius $R$ and $\eta\in C_0^\infty(B_1)$ be
such that $0\le \eta\le 1$, $\eta\equiv 1$ in $B_{R}$, $\eta\equiv 0$ outside $B_{2R}$,
and $|\nabla\eta|\le \frac{4}{R}$. For $1\le i\le k$, let
\beq{}\label{hodge_decom}
\sum_{j=1}^k P_{ji}\left\langle d\left((\phi-\bar\phi_r)\eta\right), e_j\right\rangle
=df_i+d^{\ast}g_i\  \ \hbox{ in } \mathbb R^n
\eeq
be Hodge decomposition of $\sum_{j=1}^k P_{ji}\langle d\left((\phi-\bar\phi_r)\eta\right),
e_j\rangle$ on $\mathbb R^n$,
where $f_i\in H^1\left(\mathbb R^n\right)$, $g_i\in H^1\left(\mathbb R^n,\wedge^{2}\mathbb{R}^n\right)$
is a closed $2$-form, i.e., $dg_i=0$ in $\mathbb R^n$. See Iwaniec-Martin
\cite{iwaniec-martin} for more details. Moreover, we have
the estimate
\beq{}\label{hodge_est}
\|\nabla f_i\|_{L^2(\mathbb R^n}+\|\nabla g_i\|_{L^2(\mathbb R^n)}
\le C\|d\phi\|_{L^2(B_{2R})}.\eeq
Taking $d^*$ of both sides of (\ref{hodge_decom}) and applying the equation (\ref{dirac-hm6}), we
have that for $1\le i\le k$,
\begin{eqnarray}
-\Delta f_i & =& d^*\xi_{il}\cdot\left(P_{ml}\langle d\phi, e_m\rangle\right)\ \hbox{ in }B_R, \label{fi-eqn}\\
\Delta g_i & =& dP_{ji}\wedge\left\langle d\phi,
e_j\right\rangle+P_{ji}d\phi\wedge de_j \ \ \hbox{in } B_R.
\label{gi-eqn}
\end{eqnarray}
Now we define two auxiliary $f_i^2\in H^1(B_R)$ and $g_i^2\in H^1(B_R, \wedge^2\mathbb R^n)$ on $B_R$ by
\begin{equation}\label{auxillary-fun1}
\Delta f_i^2=0 \ \hbox{ in } B_R, \  f_i^2 = f_i \ \hbox{ on } \partial B_R,
\end{equation}
\begin{equation}
\Delta g_i^2=0 \ \hbox{ in } B_R, \ g_i^2 = g_i \ \hbox{ on } \partial B_R. \label{auxillary-fun2}
\end{equation}
Set $f_i^1=f_i-f_i^2$ and $g_i^1=g_i-g_i^2$. Then $f_i^1$ and $g_i^1$ belong to $H^1_0(B_R)$.
For $1<p<\frac{n}{n-1}$, let $p'=\frac{p}{p-1}$ be its H\"older conjugate.
Recall the duality characterization of $\|\nabla u\|_{L^p(B_R)}$ for $u\in W^{1,p}_0(B_R)$:
\beq{}
\left\Vert \nabla u\right\Vert _{L^p(B_R)}\leq
C\sup \left\{\int_{B_R}\nabla u\cdot \nabla v\,dx: \
v\in W^{1,p'}_0(B_R), \  \left\Vert \nabla v\right\Vert_{L^{p'}(B_R)}\leq1\right\}.
\label{dual_lp}
\eeq
Since $p'>n$, the Sobolev embedding theorem implies
that $W_0^{1,p'}\left(B_R\right)\hookrightarrow C^{1-\frac{n}{p'}}\left(B_R\right)$
and for $v\in W_0^{1,p'}\left(B_R\right)$, with $\|\nabla v\|_{L^{p'}(B_R)}\le 1$,
there holds
\beq{}\label{sobolev_est}
\left\Vert v\right\Vert _{L^{\infty}(B_R)}\leq CR^{1-\frac{n}{p'}},
\ \left\Vert\nabla v\right\Vert_{L^2(B_R)}\le CR^{\frac{n}2-\frac{n}{p'}}.
\eeq
For any such a $v$, we can employ the equation (\ref{fi-eqn}), upon integration by parts,
use the duality between the Hardy space $\mathcal H^1$ and the BMO space
to estimate $f_i^1$, similar to Bethuel \cite{bethuel} and Rivier\'e-Struwe \cite{riviere-struwe},
as follows.
\begin{eqnarray*}
&&\int_{B_R}\nabla f_i^1\cdot\nabla v =-\int_{B_R}\Delta f_i\cdot v
=\int_{B_R} d^*\xi_{il}\cdot\left(P_{ml}\langle d\phi, e_m\rangle\right) v\\
&=&-\int_{B_R}d^*\xi_{il}\cdot d(P_{lm}e_m v)(\phi-\bar\phi_R)\\
&\le &C\|d^*\xi_{il}\cdot d(P_{lm}e_m v)\|_{\mathcal H^1(\mathbb R^n)}\left[\phi\right]_{\hbox{BMO}(B_R)}\\
&\le & C\|\nabla\xi\|_{L^2(B_R)}\left(\|\nabla P\|_{L^2(B_R)}+\|\nabla \phi\|_{L^2(B_R)}\right)
\|v\|_{L^\infty(B_R)}\left[\phi\right]_{\hbox{BMO}(B_R)}\\
&&+C\|\nabla\xi\|_{L^2(B_R)}\|\nabla v\|_{L^2(B_R)}\left[\phi\right]_{\hbox{BMO}(B_R)}\\
&\le & C\epsilon_0 R^{\frac{n-2}2}\left[R^{\frac{n-2}2}\|v\|_{L^\infty(B_R)}+\|\nabla v\|_{L^2(B_R)}\right]
\left[\phi\right]_{\hbox{BMO}(B_R)}\\
&\le &C\epsilon_0R^{\frac{n-2}{2}}[R^{1-\frac{n}{p'}+\frac{n-2}2}+R^{\frac{n}2-\frac{n}{p'}}]
\left[\phi\right]_{\hbox{BMO}(B_R)}\\
&\le & C\epsilon_0R^{\frac{n}p-1}\left[\phi\right]_{\hbox{BMO}(B_R)},
\end{eqnarray*}
where we have used that $\|\nabla e_m\|_{L^2(B_R)}\le C\|\nabla\phi\|_{L^2(B_R)}$, (\ref{small_morrey_norm}),
(\ref{morrey_est1}), and (\ref{sobolev_est}) in the derivation of these inequalities.
Taking supremum over all such $v$'s and using (\ref{dual_lp}), we obtain
\beq{}\label{fi_est1}
\left(R^{p-n}\int_{B_R}\left|\nabla f_i^1\right|^p\right)^{\frac{1}{p}}
\le C\epsilon_0 \left[\phi\right]_{\hbox{BMO}(B_R)}.
\eeq

The estimation of $g_i^1$ can be achieved in a way similar to that of $f_i^1$.
In fact, for any $v\in W^{1,p'}_0(B_R)$ satisfying (\ref{sobolev_est}), we have
\begin{eqnarray*}
&&\int_{B_R}\nabla g_i^1\cdot\nabla v =-\int_{B_R}\Delta g_i^1\cdot v=-\int_{B_R}\Delta g_i\cdot v\\
&=&-\int_{B_R} \left[dP_{ji}\wedge\left\langle d\phi,
e_j\right\rangle+P_{ji}d\phi\wedge de_j \right] v\\
&=&\int_{B_R}\left[dP_{ji}\wedge d(ve_j)+d(P_{ji}v)\wedge de_j\right]\left(\phi-\bar\phi_R\right)\\
&\le& C\left[\|dP_{ji}\wedge d(ve_j)\|_{\mathcal H^1(\mathbb R^n)}
+\|d(P_{ji}v)\wedge de_j\|_{\mathcal H^1(\mathbb R^n)}\right]\left[\phi\right]_{\hbox{BMO}(B_R)}\\
&\le&
C\|\nabla P\|_{L^2(B_R)}\left(\|\nabla v\|_{L^2(B_R)}+\|\nabla\phi\|_{L^2(B_R)}\|v\|_{L^\infty(B_R)}\right)
\left[\phi\right]_{\hbox{BMO}(B_R)}\\
&+&C\|\nabla \phi\|_{L^2(B_R)}\left(\|\nabla v\|_{L^2(B_R)}+\|\nabla P\|_{L^2(B_R)}\|v\|_{L^\infty(B_R)}\right)
\left[\phi\right]_{\hbox{BMO}(B_R)}\\
&\le & C\epsilon_0R^{\frac{n}{p}-1}\left[\phi\right]_{\hbox{BMO}(B_R)}.
\end{eqnarray*}
Taking supremum over all such $v$'s and using (\ref{dual_lp}) yields
\beq{}\label{gi_est1}
\left(R^{p-n}\int_{B_R}\left|\nabla g_i^1\right|^p\right)^{\frac{1}{p}}
\le C\epsilon_0 \left[\phi\right]_{\hbox{BMO}(B_R)}.
\eeq

Now we want to estimate $f_i^2$ and $g_i^2$. Since both $f_i^2$ and
$g_i^2$ are harmonic,  by the classical Campanato estimates for harmonic functions
(see, e.g. Giaquinta \cite{giaquinta}), (\ref{fi_est1}), and (\ref{gi_est1}),
we have that for any $0\leq r\leq R$,
it holds
\begin{eqnarray}
&&r^{p-n}\int_{B_r}\left[\left\vert \nabla f_i^2\right\vert^p
+\left\vert \nabla g_i^2\right\vert^p\right]\nonumber\\
&\le& C\left(\frac{r}{R}\right)^p\left\{R^{p-n} \int_{B_R}\left[\left\vert \nabla f_i^2\right\vert^p
+\left\vert \nabla g_i^2\right\vert^p\right]\right\}\nonumber\\
&\le &C\left(\frac{r}{R}\right)^p\left\{R^{p-n} \int_{B_R}\left[\left(\left\vert \nabla f_i\right\vert^p
+\left\vert \nabla g_i\right\vert^p\right)
+\left(|\nabla f_i^1|^p+|\nabla g_i^1|^p\right)\right]\right\}\nonumber\\
&\le & C\left(\frac{r}{R}\right)^p\left\{R^{p-n} \int_{B_R}|\nabla\phi|^p
+\epsilon_0^p \left[\phi\right]_{\hbox{BMO}(B_R)}^p\right\}. \label{fi_gi_est}
\end{eqnarray}
Therefore, using (\ref{hodge_decom}), (\ref{fi_est1}), (\ref{gi_est1}), (\ref{fi_gi_est}), and
$$|d\phi|\le \max_{i=1}^k \left|\sum_{j=1}^k P_{ji}\langle d\phi, e_j\rangle\right|,$$
we have
\begin{eqnarray}\label{phi_est}
&&r^{p-n}\int_{B_r}\left\vert \nabla \phi\right\vert^p\nonumber\\
&\le& C r^{p-n}\int_{B_r}\left[\left\vert \nabla f_i^2\right\vert^p
+\left\vert \nabla g_i^2\right\vert^p\right]+Cr^{p-n}\int_{B_r}\left[|\nabla f_i^1|^p+|\nabla g_i^2|^p\right]\nonumber\\
&\le& C\left(\frac{r}{R}\right)^p\left\{R^{p-n} \int_{B_R}|\nabla\phi|^p
+\epsilon_0^p \left[\phi\right]_{\hbox{BMO}(B_R)}^p\right\}\nonumber\\
&+&Cr^{p-n}\int_{B_R}\left[|\nabla f_i^1|^p+|\nabla g_i^2|^p\right]\nonumber\\
&\le& C\left(\frac{r}{R}\right)^p\left\{R^{p-n} \int_{B_R}|\nabla\phi|^p
+\left[1+\left(\frac{r}{R}\right)^{-n}\epsilon_0^p\right] \left[\phi\right]_{\hbox{BMO}(B_R)}^p\right\}.
\end{eqnarray}
As in \cite{riviere-struwe}, we set for $x_0\in B_1$ and $0<r\le 1$,
$$\Phi(x_0,r)=r^{p-n}\int_{B_r(x_0)}\left\vert \nabla \phi\right\vert^p,$$
and for $0<R\le 1$,
$$\Psi(R)=\sup\left\{\Phi(x_0,r):\ x_0\in B_1,\ 0<r\le R\right\}.$$
Then we have
$$\sup_{x_0\in B_1}\left[\phi\right]_{\hbox{BMO}(B_R(x_0))}^p\le C\Psi(R).$$
Thus (\ref{phi_est}) yields that there is a universal constant $C>0$ such
that for any $x_0\in B_1$ and $0<r<R\le 1$, it holds
\beq{}\label{phi_est1}
\Phi(x_0,r)\le C\left(\frac{r}{R}\right)^p
\left[1+\left(\frac{r}{R}\right)^{-n}\epsilon_0^p\right]\Psi(R).
\eeq
Now any given $\alpha\in (0,1)$, choose $\lambda\in (0,1)$ such that $2C\le\lambda^{p(\alpha-1)}$,
and choose $\epsilon_0>0$ such that $\epsilon_0^p=\lambda^n$. Then we have
\beq{} \label{phi_est2}
\Phi(x_0,\lambda R)\le 2C\lambda^p \Psi(R)\le \lambda^{p\alpha}\Psi(R)\le \lambda^{p\alpha}\Psi(R_0)
\eeq
holds for any $x_0\in B_1$, $0<R_0<1$, and $0<R\le R$. Taking supremum with respect to $x_0$ and $R<R_0$, this
gives
\beq{}\label{phi_est3}
\Psi(\lambda R_0)\le \lambda^{p\alpha}\Psi(R_0), \ \forall 0<R_0<1.
\eeq
Iteration of (\ref{phi_est3}) then yields
\beq{}\label{phi_est4}
\Psi(r)\le \left(\frac{r}{R_0}\right)^{p\alpha}\Psi(R_0), \forall 0<r\le R_0<1.
\eeq
This, combined with Morrey's decay lemma (cf. Giaquinta \cite{giaquinta}),
implies that for any $\alpha\in (0,1)$,
$\phi\in C^\alpha(B_1)$ with $[\phi]_{C^\alpha(B_1)}\le C\|\nabla\phi\|_{M^{2,2}(B_2)}$. \qed\\

Next we present a proof on the higher order regularity of $(\phi,\psi)$. The ideas are suitable
modifications of the hole-filling type argument by Giaquinta-Hildebrandt \cite{giaquinta-hildebrandt}
in the context of harmonic maps. More precisely, we have the following $C^{1,\alpha}$-regularity lemma.
\ble{}\label{higher_reg} There exist $\epsilon_0>0$ and $\theta_0\in (0,\frac12)$
such that if $(\phi,\psi)\in H^1(B_2,N)\times \mathcal S^{1,\frac43}(B_2, \mathbb C^L\tensor \phi^*TN)$
is a weakly Dirac-harmonic map satisfying
\beq{}\label{small_morrey_norm}
\|\nabla\phi\|_{M^{2,2}(B_2)}^2+\|\psi\|_{M^{4,2}(B_2)}^4\le\epsilon_0^2,
\eeq
then there exists $\mu\in (0,1)$ such that
$(\phi,\psi)\in C^{1,\mu}(B_1,N)\times C^{1,\mu}(B_1,\mathbb C^L\tensor\phi^*TN)$.
\ele
\proof The proof is divided into several steps.

\smallskip
\noindent{\it Step} 1. For any $\beta\in (0,1)$, there exists $C>0$ depending only on
$\epsilon_0$ such that
\beq{}\label{morrey_est_psi}
r^{2-n}\int_{B_r(x_0)} |\psi|^4\le Cr^{2\beta},  \ \forall \ x_0\in B_1 \ {\rm{and}}\ r\le \frac12.
\eeq
To see this, first observe that the equation $D\hspace{-0.25cm}/\hspace{0.1cm}\psi=0$ can be written as
\beq{}\label{component_dirac}
\partial\hspace{-0.2cm}/\hspace{0.1cm} \psi^i=-\Gamma_{jl}^i(\phi) \frac{\partial\phi^j}{\partial x_\alpha}
\left(\frac{\partial}{\partial x_\alpha}\circ\psi^l\right),  \ \forall 1\le i\le k,
\eeq
where $\Gamma_{jl}^i(\phi)$ is the Christoffel symbol of $(N,h)$ at $\phi$.
Note that the Lichnerowitz's formula (cf. \cite{lawson-michelsohn}) gives
$$\partial\hspace{-0.2cm}/\hspace{0.1cm}^2\psi^i=-\Delta \psi^i +\frac12 R\psi^i =-\Delta \psi^i,$$
since the scalar curvature $R=0$ on $\Omega$. Therefore, taking $\partial\hspace{-0.2cm}/\hspace{0.1cm}$
of (\ref{component_dirac})
gives
\beq{} \label{component_dirac1}
\Delta\psi^i=\partial\hspace{-0.2cm}/\hspace{0.1cm}\left[\Gamma_{jl}^i(\phi) \frac{\partial\phi^j}{\partial x_\alpha}
\left(\frac{\partial}{\partial x_\alpha}\circ\psi^l\right)\right].
\eeq
Let $\eta\in C_0^\infty(\mathbb R^n)$ be such that $0\le\eta\le 1$, $\eta=1$ in $B_{\frac{r}2}(x_0)$, and $\eta=0$ outside $B_r(x_0)$.
Define
\beq{}\label{auxi_psi}
\psi_2^i(x)=-\int_{\mathbb R^n} \frac{\partial G(x,y)}{\partial y_\beta}\frac{\partial}{\partial y_\beta}\circ
\left[\Gamma_{jl}^i(\phi) \frac{\partial\phi^j}{\partial x_\alpha}
\left(\frac{\partial}{\partial x_\alpha}\circ\psi^l\right)\eta^2\right](y)\,dy,
\eeq
where $G(x,y)=c_n |x-y|^{2-n}$ is the fundamental solution of $\Delta$ on $\mathbb R^n$. Then it is easy to see
\beq{}\label{auxi_psi1}
|\psi_2^i|(x)\le C\int_{\mathbb R^n}\frac{(|\eta\nabla\phi||\eta\psi|)(y)}{|x-y|^{n-1}}\,dy=CI_1(|\eta\nabla\phi||\eta\psi|)(x),
\eeq
where $I_1$ is the Riesz potential operator of order $1$, that is the operator whose convolution kernel is
$|x|^{1-n}, \ x\in\mathbb R^n$.
Since $\nabla\phi\in M^{2,2}(B_2)$ and $\psi\in M^{4,2}(B_2)$, by the H\"older inequality we have
$|\nabla\phi||\psi|\in M^{\frac43,2}(B_2)$. Hence
$|\eta\nabla\phi||\eta\psi|\in M^{\frac43,2}(\mathbb R^n)$  and
$$\left\||\eta\nabla\phi||\eta\psi|\right\|_{M^{\frac43,2}(\mathbb R^n)}
\le C\left\|\nabla\phi||\psi|\right\|_{M^{\frac43,2}(B_r(x_0))}
\le C\|\nabla\phi\|_{M^{2,2}(B_r(x_0))}\|\psi\|_{M^{4,2}(B_r(x_0))}.
$$
By Adams' inequality on Morrey spaces (cf. Adams \cite{adams}) $I_1: M^{\frac43,2}(\mathbb R^n)\to M^{\frac42}(\mathbb R^n)$,
we have 
\beq{}\label{psi-est1} \|\psi_2^i\|_{M^{4,2}(\mathbb R^n)}\le C\left\||\eta\nabla\phi||\eta\psi|\right\|_{M^{\frac43,2}(\mathbb R^n)}
\le C\|\nabla\phi\|_{M^{2,2}(B_r(x_0))}\|\psi\|_{M^{4,2}(B_r(x_0))}.
\eeq
By the definition of $\psi_2^i$, we have
$$\Delta\psi_2^i=\partial\hspace{-0.2cm}/\hspace{0.1cm}\left[\Gamma_{jl}^i(\phi) \frac{\partial\phi^j}{\partial x_\alpha}
\left(\frac{\partial}{\partial x_\alpha}\circ\psi^l\right)\eta^2\right]=
\partial\hspace{-0.2cm}/\hspace{0.1cm}\left[\Gamma_{jl}^i(\phi) \frac{\partial\phi^j}{\partial x_\alpha}
\left(\frac{\partial}{\partial x_\alpha}\circ\psi^l\right)\right] \ \hbox{ on }B_{\frac{r}2}(x_0)$$
so that $\Delta(\psi^i-\psi_2^i)=0$ on $B_{\frac{r}2}(x_0)$.
Hence, by the standard estimate on harmonic functions and (\ref{psi-est1}), we have
that for any $\theta\in (0,\frac12)$,
\beq{}\label{psi_est2} \left\|\psi^i-\psi_2^i\right\|_{M^{4,2}(B_{\theta r}(x_0))}\le
C\theta^{\frac12} \left\|\psi^i-\psi_2^i\right\|_{M^{4,2}(B_r(x_0))}
\le C\theta^{\frac12} \|\psi\|_{M^{4,2}(B_r(x_0))}.
\eeq
Putting (\ref{psi-est1}) and (\ref{psi_est2}) together gives
\beq{}\label{psi_est3}
\|\psi\|_{M^{4,2}(B_{\theta r}(x_0))}\le C\left(\theta^{\frac12}+\epsilon_0\right)\|\psi\|_{M^{4,2}(B_r(x_0))}.
\eeq
For any $\beta\in (0,1)$, choose $\theta=\theta(\beta)\in (0,\frac12)$ such that
$2C\le \theta^{\frac{\beta-1}2}$ and then choose $\epsilon_0$ such that
$2C\epsilon_0\le \theta^{\frac{\beta}2}$, we would have
\beq{}\label{psi-est4}
\|\psi\|_{M^{4,2}(B_{\theta r}(x_0))}\le \theta^{\frac{\beta}2}\|\psi\|_{M^{4,2}(B_r(x_0))}.
\eeq
By iteration, this clearly yields (\ref{morrey_est_psi}). \\

\noindent{\it Step} 2. For any $\beta\in (0,1)$,
any $x_0\in B_1$, and $0<r\le\frac12$, it holds
\beq{}\label{morrey_norm_decay1}
r^{2-n}\int_{B_r(x_0)}|\nabla\phi|^2\le C(\epsilon_0) r^{2\beta}.
\eeq
To do it, let $v\in H^1(B_r(x_0),\mathbb R^K)$ be such that
$$\Delta v=0 \ \hbox{ in }B_r(x_0); \ \ v= \phi \ \hbox{ on }\partial B_r(x_0).$$
Then the maximum principle and (\ref{morrey_decay1}) of Lemma \ref{morrey_decay}
imply that for any $\beta\in (0,1)$
\beq{}\label{osc_est} \left\|v-\phi\right\|_{L^\infty(B_r(x_0))}\le {\rm{osc}}_{B_r(x_0)}\phi\le Cr^\beta.
\eeq
Multiplying (\ref{dirac_harm_eqn}) by $\phi-v$ and integrating over $B_r(x_0)$ and using (\ref{osc_est}) and
(\ref{morrey_est_psi}), we have
\begin{eqnarray}\label{hole_fill1}
\int_{B_r(x_0)}\left|\nabla(\phi-v)\right|^2&\le& C\left[
\int_{B_r(x_0)}\left|\nabla\phi\right|^2\left|\phi-v\right|+\int_{B_r(x_0)}\left|\nabla\phi\right|\left|\psi\right|^2
\left|\phi-v\right|\right]\nonumber\\
&\le & Cr^{\beta}\left[\int_{B_r(x_0)}|\nabla\phi|^2+\int_{B_r(x_0)}|\psi|^4\right]\nonumber\\
&\le & Cr^\beta\int_{B_r(x_0)}|\nabla\phi|^2 +Cr^{n-2+3\beta}.
\end{eqnarray}
Hence, by the standard estimate on harmonic functions, we have
\begin{eqnarray}
(\theta r)^{2-n}\int_{B_{\theta r}(x_0)}|\nabla\phi|^2
&\le& 2 \left[(\theta r)^{2-n}\int_{B_{\theta r}(x_0)}|\nabla v|^2+(\theta r)^{2-n}\int_{B_{\theta r}(x_0)}|\nabla (\phi-v)|^2\right]\nonumber\\
&\le & C\theta^2 r^{2-n}\int_{B_r(x_0)}|\nabla\phi|^2+2(\theta r)^{2-n}\int_{B_r(x_0)}|\nabla (\phi-v)|^2\nonumber\\
&\le & C\left[\left(\theta^2+\theta^{2-n} r^{\beta}\right)r^{2-n}\int_{B_r(x_0)}|\nabla \phi|^2
+\theta^{2-n} r^{3\beta}\right]. \label{hole_fill2}
\end{eqnarray}
It is not hard to see that we can choose small $r_0>0$ and $\theta=\theta(r_0,\beta,n)\in (0,\frac12)$ such that for any $0<r\le r_0$
\beq{}\label{hole_fill3}
(\theta r)^{2-n}\int_{B_{\theta r}(x_0)}|\nabla\phi|^2
\le \theta^{2\beta}r^{2-n}\int_{B_r(x_0)}|\nabla\phi|^2+ r^{2\beta}.
\eeq
This, combined with the iteration scheme as in Giaquinta \cite{giaquinta}, implies that for any $\beta\in (0,1)$ and $x_0\in B_1$,
it holds
$$r^{2-n}\int_{B_r(x_0)}|\nabla\phi|^2\le Cr^{2\beta}, \ \forall 0<r\le r_0.$$
This yields (\ref{morrey_norm_decay1}). \\

\noindent{\it Step} 3. There exists $\mu\in (0,1)$ such that $(\phi,\psi)\in C^{1,\mu}(B_1)$. As in Step 2, let
$v\in H^1(B_r(x_0),\mathbb R^K)$ be a harmonic function with $v=\phi$ on $\partial B_r(x_0)$. Then, as in (\ref{hole_fill1}), by
using the estimates from Step 1 and Step 2 we would have that for any $\beta\in (\frac{2}{3},1)$,
\begin{eqnarray}
\label{hole_fill4}
\int_{B_r(x_0)}\left|\nabla(\phi-v)\right|^2&\le& C\left[
\int_{B_r(x_0)}\left|\nabla\phi\right|^2\left|\phi-v\right|+\int_{B_r(x_0)}\left|\nabla\phi\right|\left|\psi\right|^2
\left|\phi-v\right|\right]\nonumber\\
&\le & Cr^{\beta}\left[\int_{B_r(x_0)}|\nabla\phi|^2+\int_{B_r(x_0)}|\psi|^4\right]\nonumber\le Cr^{n-2+3\beta}.
\end{eqnarray}
Hence, using the Campanato estimate for harmonic functions, we have
\begin{eqnarray}
&&(\theta r)^{-n}\int_{B_{\theta r}(x_0)}\left|\nabla \phi-\overline{\nabla\phi}_{x_0,\theta r}\right|^2\nonumber\\
&\le & 2 \left[(\theta r)^{-n}\int_{B_{\theta r}(x_0)}\left|\nabla v-\overline{\nabla v}_{x_0,r}\right|^2
+(\theta r)^{-n}\int_{B_{\theta r}(x_0)}\left|\nabla (\phi-v)\right|^2\right]\nonumber\\
&\le & C\left[\theta^2 r^{-n}\int_{B_r(x_0)}\left|\nabla \phi-\overline{\nabla\phi}_{x_0,r}\right|^2+\theta^{-n}r^{3\beta-2}\right]\nonumber\\
&\le &\frac{1}{2}r^{-n}\int_{B_r(x_0)}\left|\nabla \phi-\overline{\nabla\phi}_{x_0,r}\right|^2+Cr^{2\mu}
\end{eqnarray}
provide that we first choose $\theta\in (0,\frac12)$ sufficiently small, and then choose $r_0$ so small that
$C\theta^{-n}r^{3\beta-2}\le r^{2\mu}$ for $0<r\le r_0$,
where $$\overline{\nabla\phi}_{x_0, r}=\frac{1}{|B_r(x_0)|}\int_{B_r(x_0)}\nabla \phi$$ is the average of $\nabla \phi$ over
$B_r(x_0)$. It follows from the same iteration scheme as in \cite{giaquinta} that
$$r^{-n}\int_{B_r(x_0)}|\nabla\phi-\overline{\nabla\phi}_{x_0,r}|^2
\le Cr^{2\mu}, \ \forall x_0\in B_1, \ 0<r\le r_0.$$
This, combined with the characterization of $C^\mu$ by the Campanato space, implies $\nabla\phi\in C^{\mu}(B_1)$.
Substituting $\nabla\phi\in C^\mu(B_1)$ into the equation (\ref{component_dirac}), one can easily conclude that
$\psi\in C^{1,\mu}(B_1)$.
\qed\\

At the end of this section, we complete the proof of Theorem \ref{small_reg}.

\medskip
\noindent{\bf Completion of Proof of Theorem \ref{small_reg}}.

\smallskip
Combining Lemma \ref{morrey_decay} and Lemma \ref{higher_reg}, we know that $(\phi,\psi)\in C^{1,\mu}(B_{\frac{r_0}2}(x_0))$.
The higher order regularity then follows from the standard bootstrap argument for both equations (\ref{dirac_harm_eqn})
and (\ref{dirac_eqn}). We omit the details. \qed

\section{Stationary Dirac-harmonic maps}

In this section, we introduce the notion of stationary Dirac-harmonic maps, which is a natural extension
of stationary harmonic maps. Any smooth Dirac-harmonic map is a stationary Dirac-harmonic
map (see \cite{chen-jost-li-wang1}). We prove several interesting properties, including a partial regularity theorem, for
stationary Dirac-harmonic maps. To simplify the presentation, we assume throughout this section
that $(M,g)=(\Omega,g_0)$, where $\Omega\subseteq\mathbb R^n$ is a bounded smooth domain and
$g_0$ is the Euclidean metric on $\mathbb R^n$. Thus the spinor bundle $\Sigma$ associated with
$M$ can also be identified by $\Sigma=\Omega\times \mathbb C^L$, $L={\rm{rank}}_{\mathbb C}\Sigma$. We remark that
one can modify the proofs of Lemma \ref{free-stress}, Lemma \ref{mono_id1}, and Proposition \ref{mono_id5} without
difficulties so that (\ref{stationary_id1}) and (\ref{mono_id012}) in \S1 hold for any general Riemannian manifold $(M,g)$,
we leave the details to the interested readers.

\bde{} A weakly Dirac-harmonic map $(\phi, \psi)\in H^1(\Omega,N)\times {\mathcal S}^{1,\frac{4}{3}}
(\Omega, \mathbb C^L\tensor \phi^*TN))$
is called to be a stationary Dirac-harmonic map, if it is also a critical point of $L(\phi,\psi)$ with respect to the
domain variations, i.e., for any $Y\in C^\infty_0(\Omega,\mathbb R^n)$, it holds
\beq{} \label{stationary-dirac}
\frac{d}{dt}\big|_{t=0}\left[\int_\Omega \left(|\nabla\phi_t|^2+\langle \psi_t, D\hspace{-0.25cm}/\hspace{0.1cm}\psi_t\rangle\right)\right]=0,
\eeq
where $\phi_t(x)=\phi(x+tY(x))$ and $\psi_t(x)=\psi(x+tY(x))$.
\ede

We now derive the stationarity identity for stationary Dirac-harmonic maps defined as above.
\ble{}\label{free-stress} Let $(\phi, \psi)\in H^1(\Omega,N)\times {\mathcal S}^{1,\frac43}(\Omega, \mathbb C^L\tensor \phi^*TN))$
be a weakly Dirac-harmonic map. Then  $(\phi,\psi)$ is a stationary Dirac-harmonic map iff
for any $Y\in C_0^\infty(\Omega,\mathbb R^n)$, it holds
\beq{}\label{stat_id1}
\int_\Omega \left[\langle \frac{\partial\phi}{\partial x_\alpha}, \frac{\partial\phi}{\partial x_\beta}\rangle
-\frac12|\nabla \phi|^2\delta_{\alpha\beta}+\frac12{\rm{Re}}\left\langle\psi, \frac{\partial}{\partial x_\alpha}
\circ \nabla_{\frac{\partial}{\partial x_\beta}}\psi\right\rangle\right]
\frac{\partial Y_\alpha}{\partial x_\beta}=0.
\eeq
\ele
\proof For $t\in\mathbb R$ with small $|t|$,
denote $y=F_t(x)=x+tY(x):\Omega\to\Omega$ and $x=F_t^{-1}(y)$.
It is a standard calculation (see, e.g. \cite{lin-wang}) that
\beq{} \label{stat_id2}
\frac{d}{dt}\big|_{t=0}\int_\Omega \left|d\phi_t\right|^2\,dx
=\int_\Omega \left[2\langle \frac{\partial\phi}{\partial x_\alpha}, \frac{\partial\phi}{\partial x_\beta}\rangle
\frac{\partial Y_\alpha}{\partial x_\beta}-|\nabla\phi|^2\hbox{div}(Y)\right].
\eeq
Now we compute $\frac{d}{dt}\big|_{t=0}\int_\Omega  \langle \psi_t, D\hspace{-0.25cm}/\hspace{0.1cm}\psi_t\rangle$.
First, by remark 1.3, we have
$$
\int_\Omega  \left\langle \psi_t, D\hspace{-0.25cm}/\hspace{0.1cm}\psi_t\right\rangle
=\int_\Omega  {\rm{Re}}\left\langle \psi_t,D\hspace{-0.25cm}/\hspace{0.1cm}\psi_t\right\rangle.$$
Before taking $\frac{d}{dt}$, we perform a change of variable as follows.
For $y=F_t(x)$, since
$$\frac{\partial}{\partial x_\alpha}=\frac{\partial y_\beta}{\partial x_\alpha}\frac{\partial}{\partial y_\beta},$$
we have
$$D\hspace{-0.25cm}/\hspace{0.1cm}\psi_t=\frac{\partial}{\partial x_\alpha}(F_t(x))\circ
\nabla_{\frac{\partial}{\partial x_\alpha}(F_t(x))}\psi
=\frac{\partial y_\beta}{\partial x_\alpha}
\frac{\partial}{\partial x_\alpha}(y)\circ\nabla_{\frac{\partial}{\partial y_\beta}}\psi.$$
Thus
$$
\int_\Omega  {\rm{Re}}\left\langle \psi_t, D\hspace{-0.25cm}/\hspace{0.1cm}\psi_t\right\rangle
=\int_\Omega \sum_{\alpha,\beta}
\frac{\partial y_\beta}{\partial x_\alpha} {\rm{Re}}\left\langle \psi, \frac{\partial}{\partial x_\alpha}(y)\circ
\nabla_{\frac{\partial}{\partial y_\beta}}\psi\right\rangle \hbox{Jac}F_t^{-1}\,dy.
$$
Since
$$\frac{d}{dt}\big|_{t=0}\hbox{Jac}F_t^{-1}=-\hbox{div}(Y)
\ {\rm{and}} \ \ \frac{d}{dt}\big|_{t=0}\frac{\partial y_\beta}{\partial x_\alpha}
=\frac{\partial Y_\beta}{\partial x_\alpha},$$
we have
\begin{eqnarray}
&&\frac{d}{dt}\big|_{t=0}\int_\Omega {\rm{Re}}\left\langle \psi_t, D\hspace{-0.25cm}/\hspace{0.1cm}\psi_t\right\rangle\nonumber\\
&=&\int_\Omega{\rm{Re}}\left\langle\psi,D\hspace{-0.25cm}/\hspace{0.1cm}\psi\right\rangle\left[\frac{d}{dt}|_{t=0}\hbox{Jac}(F_t^{-1})\right]
+\int_\Omega \sum_{\alpha,\beta}{\rm{Re}}\left\langle\psi,
\frac{\partial}{\partial x_\alpha}\circ\nabla_{\frac{\partial}{\partial x_\beta}}\psi\right\rangle
\frac{\partial Y_\beta}{\partial x_\alpha}\nonumber\\
&=&-\int_\Omega {\rm{Re}}\left\langle\psi,D\hspace{-0.25cm}/\hspace{0.1cm}\psi\right\rangle \hbox{div}(Y)
+\int_\Omega \sum_{\alpha,\beta} {\rm{Re}}
\left\langle \psi, \frac{\partial}{\partial x_\alpha}\circ \nabla_{\frac{\partial}{\partial x_\beta}}\psi\right\rangle
\frac{\partial Y_\beta}{\partial x_\alpha}\nonumber\\
&=& \int_\Omega {\rm{Re}}\left\langle \psi, \frac{\partial}{\partial x_\alpha}
\circ \nabla_{\frac{\partial}{\partial x_\beta}}\psi\right\rangle
\frac{\partial Y_\beta}{\partial x_\alpha}, \label{stat_id3}
\end{eqnarray}
where we have used $D\hspace{-0.25cm}/\hspace{0.1cm}\psi=0$ in the last step. It is clear that (\ref{stat_id1})
follows from (\ref{stat_id2}) and (\ref{stat_id3}).
\qed\\

It is well known that any stationary harmonic map $u:\mathbb R^n\to N$ with finite Dirichlet energy is
a constant map. Here we prove that the same conclusion holds for stationary Dirac-harmonic maps.

\bth{} \label{liouville} For $n\ge 3$, assume that $(\phi,\psi)\in H^1_{\hbox{loc}}(\mathbb R^n, N)
\times \mathcal S^{1,\frac43}_{\hbox{loc}}\left(\Gamma(\Sigma\times\phi^*TN)\right)$
is a stationary Dirac-harmonic map satisfying
\beq{}\int_{\mathbb R^n}\left(|\nabla\phi|^2+|\nabla\psi|^{\frac43}+|\psi|^4\right)<+\infty.
\label{finite}
\eeq
Then $\phi\equiv\hbox{constant}$ and $\psi\equiv 0.$
\eth
\proof  For any large $R>0$, let $\eta\in C^\infty_0(\mathbb R)$ be such that $\eta\equiv 1$ for $|t|\le R$,
$\eta\equiv 0$ for $|t|\ge 2R$, and $|\eta'(t)|\le \frac{4}{R}$.
Let $Y(x)=x\eta(|x|)\in C^\infty_0(\mathbb R^n,\mathbb R^n)$. Note that
$$\frac{\partial Y_\beta}{\partial x_\alpha}=\eta(|x|)\delta_{\alpha\beta}+\frac{x_\alpha x_\beta}{|x|}\eta'(|x|).$$
Substituting $Y$ into the stationarity identity (\ref{stat_id1}) yields
\begin{eqnarray*}
&&\left(\frac{n}2-1\right)\int_{\mathbb R^n} |\nabla\phi(x)|^2\eta(|x|)-\frac12\int_{\mathbb R^n}{\rm {Re}}\left\langle\psi,
\frac{\partial}{\partial x_\alpha}\circ\nabla_{\frac{\partial}{\partial x_\alpha}}\psi\right\rangle \eta(|x|) \\
&=&\int_{\mathbb R^n}\left[\langle \frac{\partial\phi}{\partial x_\alpha},
\frac{\partial\phi}{\partial x_\beta}\rangle-\frac12|\nabla\phi|^2\delta_{\alpha\beta}
+\frac12{\rm{Re}}\left\langle \psi,
\frac{\partial}{\partial x_\alpha}\circ\nabla_{\frac{\partial}{\partial x_\beta}}\psi\right\rangle\right]
\eta'(|x|)\frac{x_\alpha x_\beta}{|x|}\\
&=&\int_{\mathbb R^n}\left[\left(|\frac{\partial \phi}{\partial r}|^2-\frac12|\nabla \phi|^2\right)
+\frac12{\rm{Re}}\left\langle\psi,\frac{\partial}{\partial r}\circ\nabla_{\frac{\partial}{\partial r}}\psi\right\rangle\right]
|x|\eta'(|x|). \end{eqnarray*}
Since $D\hspace{-0.25cm}/\hspace{0.1cm}\psi=0$, we have
$$\frac12\int_{\mathbb R^n}{\rm {Re}}\left\langle\psi,
\frac{\partial}{\partial x_\alpha}\circ\nabla_{\frac{\partial}{\partial x_\alpha}}\psi\right\rangle \eta(|x|)
=\frac12\int_{\mathbb R^n}{\rm {Re}}\left\langle\psi,
D\hspace{-0.25cm}/\hspace{0.1cm}\psi\right\rangle \eta(|x|)=0.$$
The right hand side can be estimated by
\begin{eqnarray}
&&\left|\int_{\mathbb R^n}\left[\left(|\frac{\partial \phi}{\partial r}|^2-\frac12|\nabla \phi|^2\right)
+\frac12{\rm{Re}}\left\langle\psi,\frac{\partial}{\partial r}\circ\nabla_{\frac{\partial}{\partial r}}\psi\right\rangle\right]
|x|\eta'(|x|)\right|\nonumber\\
&\le& C\int_{B_{2R}\setminus B_R}\left[|\nabla\phi|^2+|\psi||\nabla\psi|\right]\nonumber\\
&\le & C\int_{B_{2R}\setminus B_R}\left[|\nabla\phi|^2+|\psi|^4+|\nabla\psi|^{\frac43}\right]. \label{rhs}
\end{eqnarray}
On the other hand, the left hand side is bounded below by
\beq{}\left(\frac{n}2-1\right)\int_{\mathbb R^n} |\nabla\phi(x)|^2\eta(|x|)
\ge \left(\frac{n}{2}-1\right)\int_{B_R}|\nabla\phi|^2.\label{lhs}
\eeq
Since $n\ge 3$, (\ref{rhs}) and (\ref{lhs}) imply
\beq{}\label{lhs_rhs}
\int_{B_R}|\nabla\phi|^2\le C\int_{B_{2R}\setminus B_R}\left[|\nabla\phi|^2+|\psi|^4+|\nabla\psi|^{\frac43}\right],
\eeq
this and (\ref{finite}) imply, after sending $R$ to $\infty$,
$$\int_{\mathbb R^n}|\nabla\phi|^2=0, $$
i.e., $\phi\equiv \hbox{constant}$. Substituting this $\phi$ into the equation of $\psi$, we have
$\partial\hspace{-0.2cm}/\hspace{0.1cm}\psi=0\ \hbox{ in }\mathbb R^n.$
Since $\psi\in L^4(\mathbb R^n)$, it follows easily that $\psi\equiv 0$. \hfill\qed\\

Now we derive an identity for stationary Dirac-harmonic maps, that is similar to the
monotonicity identity for stationary harmonic maps.

\ble{}\label{mono_id1} Assume that $\left(\phi,\psi\right)\in H^1(\Omega,N)\times \mathcal S^{1,\frac43}(\Omega,
\mathbb C^L\tensor \phi^*TN)$ is a stationary Dirac-harmonic map.
Then for any $x_0\in \Omega$ and $0<R_1\le R_2<{\rm{dist}}(x_0,\partial\Omega)$, it holds
\begin{eqnarray}
&&R_2^{2-n}\int_{B_{R_2}(x_0)}\left|\nabla\phi\right|^2\,dx
-R_1^{2-n}\int_{B_{R_1}(x_0)}\left\vert \nabla \phi\right\vert^{2}\,dx\nonumber\\
&=&\int_{R_1}^{R_2}r^{2-n}\left(\int_{\partial B_r(x_0)}\left[ 2\left\vert \frac{%
\partial \phi}{\partial r}\right\vert ^{2}+{\rm{Re}}\left\langle \psi ,\frac{\partial }{%
\partial r}\circ \nabla _{\frac{\partial }{\partial r}}\psi \right\rangle\right] dH^{n-1}\right)\,dr,
\label{mono_id2}
\end{eqnarray}
where $\frac{\partial}{\partial r}=\frac{\partial}{\partial |x-x_0|}$.
\end{lemma}
\proof For simplicity, assume $x_0=0\in\Omega$. The argument is similar to that of Theorem
\ref{liouville}. For completeness, we outline it again. For $\epsilon>0$ and $0<r<{\rm{dist}}(0,\partial\Omega)$,
let $\eta_\epsilon(x)=
\eta _{\epsilon }(|x|) \in C_{0}^{\infty }\left( B_{r}\right)$ be such that
$0\le\eta_\epsilon\le 1$, $\eta_{\epsilon}=1$ for $\left\vert x\right\vert \leq r\left( 1-\epsilon \right)$.
Choose $Y(x)=x\eta_{\varepsilon}\left(\left\vert x\right\vert\right)$.
Note
$$\frac{\partial Y_\beta}{\partial x_\alpha}=\delta_{\alpha\beta}\eta_\epsilon(|x|)
+\frac{x_\alpha x_\beta}{|x|}\eta_\epsilon'(|x|).$$
Substituting $Y$ into (\ref{stat_id1}), we have
\begin{eqnarray*}
&&\left(2-n\right)\int_{B_r} |\nabla\phi(x)|^2\eta_\epsilon(|x|)+
\int_{B_r}{\rm {Re}}\left\langle\psi, \frac{\partial}{\partial x_\alpha}\circ\nabla_{\frac{\partial}{\partial x_\alpha}}\psi\right\rangle \eta_\epsilon(|x|) \\
&=&-\int_{B_r}\left[2\langle \frac{\partial\phi}{\partial x_\alpha},
\frac{\partial\phi}{\partial x_\beta}\rangle-|\nabla\phi|^2\delta_{\alpha\beta}
+{\rm{Re}}\left\langle \psi,
\frac{\partial}{\partial x_\alpha}\circ\nabla_{\frac{\partial}{\partial x_\beta}}\psi\right\rangle\right]
\eta_\epsilon'(|x|)\frac{x_\alpha x_\beta}{|x|}\\
&=&-\int_{B_r}\left[\left(2|\frac{\partial \phi}{\partial r}|^2-|\nabla \phi|^2\right)
+{\rm{Re}}\left\langle\psi,\frac{\partial}{\partial r}\circ\nabla_{\frac{\partial}{\partial r}}\psi\right\rangle\right]
|x|\eta_\epsilon'(|x|).
\end{eqnarray*}
Using the equation $D\hspace{-0.25cm}/\hspace{0.1cm}\psi=\frac{\partial}{\partial x_\alpha}\circ\nabla_{\frac{\partial}{\partial x_\alpha}}\psi=0$
and sending $\epsilon$ to $0$, this yields
\begin{eqnarray}
&&\left(2-n\right)\int_{B_r}\left\vert \nabla \phi\right\vert^2\,dx
+r\int_{\partial B_{r}}\left\vert \nabla \phi\right\vert^2\,dH^{n-1}\notag \\
&=&2r\int_{\partial B_{r}}\left\vert \frac{\partial \phi}{\partial r}\right\vert^{2}
dH^{n-1}+r\int_{\partial B_{r}}{\rm{Re}}\left\langle \psi ,
\frac{\partial}{\partial r}\circ \nabla _{\frac{\partial }{\partial r}}\psi \right\rangle \,dH^{n-1},
\label{mono_id3}
\end{eqnarray}
or equivalently,
$$
\frac{d}{dr}\left( r^{2-n}\int_{B_{r}}\left\vert \nabla \phi\right\vert
^{2}\,dx\right) =r^{2-n}\int_{\partial B_{r}}
\left[ 2\left\vert \frac{\partial \phi}{\partial r}\right\vert ^{2}
+{\rm{Re}}\left\langle \psi ,\frac{\partial }{\partial r}\circ \nabla _{\frac{\partial }{\partial r}}\psi
\right\rangle\right] \,dH^{n-1}.
$$
Integrating $r$ from $R_1$ to $R_2$ yields (\ref{mono_id1}). \qed\\

In contrast with stationary harmonic maps, (\ref{mono_id2}) doesn't imply that the renormalized
energy $$R^{2-n}\int_{B_R(x_0)}|\nabla\phi|^2$$ is monotone increasing with respect to $R$ yet.
In order to have such a monotonicity property, we need to assume that $\nabla\psi$ has
higher integrability. More precisely, we have
\bpr{} \label{mono_id5}
Assume that $\left(\phi,\psi\right)\in H^1(\Omega,N)\times \mathcal S^{1,\frac43}(\Omega,
\mathbb C^L\tensor \phi^*TN)$ is a stationary Dirac-harmonic map. If, in additions,
$\nabla\psi\in L^p(\Omega)$ for some $\frac{2n}{3}<p\le n$, then there exists $C_0>0$
depending only on $\|\nabla\psi\|_{L^p(\Omega)}$ such that for any $x_0\in\Omega$
and $0<R_1<R_2<{\rm{dist}}(x_0,\partial\Omega)$, it holds
\beq{}\label{mono_id6}
R_1^{2-n}\int_{B_{R_1}(x_0)}|\nabla\phi|^2\le R_2^{2-n}\int_{B_{R_2}(x_0)}|\nabla\phi|^2+
C_0 R_2^{3-\frac{2n}{p}}.
\eeq
\epr
\proof For simplicity, assume $x_0=0$.
For $x\in \Omega$, denote
$$f(x)=\left|{\rm{Re}}\left\langle \psi ,\frac{\partial }
{\partial r}\circ \nabla _{\frac{\partial }{\partial r}}\psi \right\rangle(x)\right|.$$
Since $\nabla\psi\in L^p(\Omega)$, by Sobolev embedding theorem
we have $\psi\in L^{\frac{np}{n-p}}(\Omega)$. Since $f(x)\le C|\psi||\nabla\psi|$,
the H\"older inequality implies that $f\in L^{q}(\Omega)$ with $q=\frac{np}{2n-p}$.
Since $p>\frac{2n}3$, it is easy to see that $q>\frac{n}2$. It then follows that
for any $R<R_0={\rm{dist}}(0,\partial\Omega)$,
\begin{eqnarray}\label{lq_est1}\int_0^{R}r^{1-n}\int_{B_r}f(x)\,dx
&\le& \left(\int_0^{R}r^{1-\frac{n}{q}}\,dr\right) \|f\|_{L^q(B_{R_0})}\nonumber\\
&=&\left(\frac{q}{2q-n}\right)R^{2-\frac{n}{q}}\|f\|_{L^q(B_{R_0})}<+\infty,
\end{eqnarray}
and
\beq{}\label{lq_est2}
R^{2-n}\int_{B_{R}}f\le  R^{2-\frac{n}{q}}\|f\|_{L^q(B_{R_0})}.
\eeq
For any $0<R_1\le r\le R_2<R_0$, set
$$g(r)=\int_{\partial B_r}f(x)\,dH^{n-1}.$$
Then, by integration by parts, we have
\begin{eqnarray*}
&&\int_{R_1}^{R_2} r^{2-n}\int_{\partial B_r}\left|{\rm{Re}}\left\langle \psi ,\frac{\partial }
{\partial r}\circ \nabla _{\frac{\partial }{\partial r}}\psi \right\rangle\right|\,dH^{n-1}\\
&=&\int_{R_1}^{R_2}r^{2-n}g(r)\,dr=\int_{R_1}^{R_2}r^{2-n}\,d\left(\int_{B_r}f(x)\,dx\right)\\
&=& R_2^{2-n}\int_{B_{R_2}}f(x)\,dx-R_1^{2-n}\int_{B_{R_1}}f(x)\,dx
+(n-2)\int_{R_1}^{R_2} r^{1-n}\left(\int_{B_r}f(x)\,dx\right).
\end{eqnarray*}
This, combined with (\ref{mono_id2}), then implies
\begin{eqnarray*}
&&R_2^{2-n}\int_{B_{R_2}}|\nabla\phi|^2+R_2^{2-n}\int_{B_{R_2}}f
+(n-2)\int_0^{R_2}r^{1-n}\int_{B_r} f\\
&&\ge R_1^{2-n}\int_{B_{R_1}}|\nabla\phi|^2+R_1^{2-n}\int_{B_{R_1}}f
+(n-2)\int_0^{R_1}r^{1-n}\int_{B_r} f\\
&&+\int_{R_1}^{R_2}r^{2-n}\int_{\partial B_r}\left|\frac{\partial\phi}{\partial r}\right|^2\,dH^{n-1}.
\end{eqnarray*}
It is easy to see that this inequality, (\ref{lq_est1}), and (\ref{lq_est2})
imply (\ref{mono_id6}). \qed\\

With the help of Proposition \ref{mono_id5} and Theorem \ref{small_reg}, we can prove
Theorem \ref{nd_reg}. \\

\noindent{\bf Proof of Theorem \ref{nd_reg}}.

For simplicity, assume $M=\Omega\subseteq\mathbb R^n$ and $g=g_0$ is the Euclidean metric
on $\mathbb R^n$. Since $\nabla\psi\in L^p(\Omega)$ for some
$p>\frac{2n}3$, we have by Sobolev's embedding theorem that
$\psi\in L^{q}(\Omega)$ for $q=\frac{np}{n-p}>2n$. Hence, for any ball $B_R(x)\subseteq\Omega$,
by the H\"older inequality we have
\beq{}\label{l4_est}
R^{2-n}\int_{B_R(x)}|\psi|^4\le \left(\int_{B_R(x)}|\psi|^q\right)^{\frac{4}{q}} R^{2-\frac{4n}q}
\le \|\nabla\psi\|^4_{L^p(\Omega)}R^{2-\frac{4n}q}.
\eeq

Let $\epsilon_0>0$ be the constant given by Theorem \ref{small_reg}. For a large constant $C(n)>0$ to be chosen
later, define
$$\mathcal S(\phi)
=\bigcap_{R>0}\left\{x\in\Omega: R^{2-n}\int_{B_R(x)}|\nabla\phi|^2>\frac{\epsilon_0^2}{C(n)}\right\}.$$
It is well known (cf. Evans-Gariepy \cite{evans-gariepy}) that $H^{n-2}\left(\mathcal S(\phi)\right)=0$.
For any $x_0\in\Omega\setminus\mathcal S(\phi)$, there exists $r_0>0$ such that
$$\left(2r_0\right)^{2-n}\int_{B_{2r_0}(x_0)}\left|\nabla\phi\right|^2\le\frac{\epsilon_0^2}{C(n)}.$$
Hence
$$\sup_{x\in B_{r_0}(x_0)}\left\{r_0^{2-n}\int_{B_{r_0}(x)}\left|\nabla\phi\right|^2\right\}
\le \frac{2^{n-2}\epsilon_0^2}{C(n)}.$$
Applying the monotonicity inequality (\ref{mono_id6}), this implies
\beq{}\label{morrey_small1}\sup_{x\in B_{r_0}(x_0), 0<r\le r_0}
\left\{r^{2-n}\int_{B_{r}(x)}\left|\nabla\phi\right|^2\right\}
\le \frac{2^{n-2}\epsilon_0^2}{C(n)}+C_0 r_0^{3-\frac{2n}{p}}
\le\frac{\epsilon_0^2}{4},\eeq
provide that we choose $C(n)>2^{n+1}$ and $r_0\le \left(\frac{\epsilon_0^2}{8C_0}\right)^{\frac{p}{3p-2n}}$.
On the other hand, by (\ref{l4_est}) we have
\beq{}\label{morrey_small2} \sup_{x\in B_{r_0}(x_0), 0<r\le r_0}
\left\{r^{2-n}\int_{B_{r}(x)}\left|\psi\right|^4\right\}\le \|\nabla\psi\|_{L^p(\Omega)}^4 r_0^{\frac{6p-4n}{p}}
\le \frac{\epsilon_0^2}{4},\eeq
provide that we choose $r_0<\left(\frac{\epsilon_0^2}{4\|\nabla\psi\|_{L^p}^4}\right)^{\frac{p}{6p-4n}}$.
Combining (\ref{morrey_small1}) with (\ref{morrey_small2}), we have that there exists $r_0>0$ sufficiently small
such that
\beq{}\label{morrey_small3}
\sup_{x\in B_{r_0}(x_0), 0<r\le r_0}
\left\{r^{2-n}\int_{B_{r}(x)}\left|\nabla\phi\right|^2+r^{2-n}\int_{B_r(x)}|\psi|^4\right\}
\le\frac{\epsilon_0^2}{2}.
\eeq
Thus Theorem \ref{small_reg} implies that $(\phi,\psi)\in C^\infty(B_{\frac{r_0}2}(x_0),N)
\times C^\infty(B_{\frac{r_0}2}(x_0), \mathbb C^L\tensor\phi^*TN)$. Note that this also yields $\Omega\setminus\mathcal S(\phi)$
is an open set. The proof of Theorem \ref{nd_reg} is now complete.
\qed\\

\section{Convergence of approximate Dirac-harmonic maps}

In dimension two, the weak convergence theorem of approximated harmonic maps or Palais-Smale sequences of
Dirichlet energy functional for maps into Riemannian manifolds was first proved by Bethuel \cite{bethuel1}.
Subsequently, alternative proofs were given by Freire-M\"uller-Struwe \cite{FMS}, and Wang \cite{wang} by
employing the moving frame and various techniques including the concentration compactness method. Very recently,
Rivier\'e \cite{riviere} gave another proof using the conservation laws.

In this section, we extend such a convergence theorem to sequences of approximate Dirac-harmonic maps
from a spin Riemann surface. The key ingredient is to first use the moving frame to rewrite the equation of
approximate harmonic maps into the form similar to (\ref{new_dirac_harm_eqn}), and then use Rivier\'e's Coulomba
gauge construction technique to rewrite it into the form (\ref{dirac-hm5})
in which the concentration compactness method similar to that of \cite{FMS} to pass to the limit.

Since $\int |d\phi|^2+|\psi|^4$ is conformally invariant in dimension two, it follows from a scaling argument
and a covering argument that Theorem \ref{conv_dirac_hm} follows from the following lemma. For simplicity,
we assume that $(M,g)=(\Omega,g_0)$ for some bounded smooth domain $\Omega\subseteq\mathbb R^2$ with
the Euclidean metric $g_0$. Denote by $B_r\subseteq\mathbb R^2$ be the ball center at $0$ with radius $r$.

\ble{} \label{small_convergence} There exists $\epsilon_1>0$ such that if $(\phi_p,\psi_p)\in H^1(B_1, N)\times
\mathcal S^{1,\frac43}(B_1, \mathbb C^2\tensor\phi_m^*TN)$ is a sequence of approximate Dirac-harmonic maps, i.e.
$$\tau(\phi_p)=\mathcal R^N(\phi_p,\psi_p)+u_p; \ D\hspace{-0.25cm}/\hspace{0.1cm}\psi_p= v_p \ \hbox{ on }B_1,$$
$$u_p\rightarrow 0 \hbox{ strongly in }H^{-1}(B_1) \ \hbox{ and } v_p \rightharpoonup 0 \hbox{ in }L^{\frac43}(B_1),$$
and
\beq{}\label{small_energy_assump}
\int_{B_1}\left(|\nabla\phi_p|^2+|\psi_p|^4\right)\le\epsilon_1^2.
\eeq
If $(\phi_p,\psi_p)\rightharpoonup (\phi,\psi) \hbox{ in } H^1(B_1,N)\times
\mathcal S^{1,\frac43}(B_1, \mathbb C^2\tensor\mathbb R^K)$, then $(\phi,\psi)
\in H^1(B_1,N)\times \mathcal S^{1,\frac43}(B_1, \mathbb C^2\tensor\phi^*TN)$
is a weakly Dirac harmonic map.
\ele
\proof First observe that the argument of Proposition \ref{invariance} can easily be modified to show
that if $(\widetilde N,\widetilde h)$ is another Riemannian manifold
and $f:(N,h)\to (\widetilde N,\widetilde h)$ is a totally geodesic, isometric embedding,
and if we set $\widetilde \phi=f(\phi)$ and $\widetilde\psi=f_*(\psi)$, then
$(\widetilde \phi, \widetilde\psi)\in H^1(B_1,\widetilde N)\times
\mathcal S^{1,\frac43}(B_1,\mathbb C^2\tensor (\widetilde\phi)^*T\widetilde N)$ is a
sequence of approximate harmonic maps with $(u_p,v_p)$ replaced by
$(\widetilde u_p,\widetilde v_p)$, where $\widetilde u_p=f_*(u_p)$ and
$\widetilde v_p=f_*(v_p)$. Moreover, it is easy to check that
$$\widetilde u_p\rightarrow 0 \hbox{ strongly in }H^{-1}(B_1)
\ \hbox{ and } \widetilde v_p\rightharpoonup 0 \ \hbox{ in } L^{\frac43}(B_1),$$
$$\widetilde \phi_p\rightharpoonup \widetilde \phi=f(\phi) \hbox{ in }H^1(B_1);
\ \widetilde \psi_p \rightharpoonup \widetilde \psi=f_*(\psi) \hbox{ in }\mathcal S^{1,\frac43}(B_1),$$
and
$$\int_{B_1}|\nabla \widetilde \phi_p|^2+|\widetilde\psi_p|^2\le \epsilon_0^2.$$
With this reduction, we may assume that there exists global orthonormal frame $\{\hat e_i\}_{i=1}^k$ on $(N,h)$.
For any $p$, let $e_i^p=\hat e_i(\phi_p)$, $1\le i\le k$, be the orthonormal frame along with $\phi_p$.
Then, similar to Lemma \ref{dirac_hm_eqn}, Proposition \ref{skew-sym}, and (\ref{dirac-hm4}),  we have
\beq{}\label{approx_dirac_hm}
d^*\left(\langle d\phi_p,e_i^p\rangle\right)=\sum_{j=1}^k \Theta_{ij}^p \langle d\phi_p,e_j^p\rangle+u_p^i,
\eeq
where
\beq{}\label{skew-form1}
\Theta_{ij}^p=\Omega_{ij}^p+\langle de_i^p, e_j^p\rangle, \ 1\le i, j\le k;\ \  u_p^i=\langle u_p, e_i^p\rangle,
\ 1\le i\le k,\eeq
and
\beq{}\label{skew-form2}
\Omega_{ij}^p=\sum_{\alpha=1}^n \left[
\sum_{l,m=1}^k R^N(\phi_p)\left(e_i^p, e_j^p, e_l^p, e_m^p\right)\langle \psi_p^m,
\frac{\partial}{\partial x_\alpha}\circ \psi^l_p\rangle\right]\,dx_\alpha, 1\le i, j\le k.
\eeq
Since $\Theta^p=\left(\Theta_{ij}^p\right)$ satisfies
$$\int_{B_1}\left|\Theta^p\right|^2\le C\int_{B_1}\left(|\nabla\phi_p|^2+|\psi_p|^4\right)\le C\epsilon_1^2
\le\epsilon_0^2,$$
provide $\epsilon_1\le \frac{\epsilon_0}{\sqrt C}$, where $\epsilon_0$ is the same constant as
in Lemma \ref{coulomb_gauge}. Hence Lemma \ref{coulomb_gauge} implies that there
exist $Q^p\in H^1(B_1, {\rm{SO}}(k))$ and $\xi^p\in H^1(B_1, {\rm{so}}(k)\tensor\wedge^2\mathbb R^2)$ such that
\beq{}\label{coulomb_gauge_p1}
\left(Q^p\right)^{-1}dQ^p-\left(Q^p\right)^{-1} \Theta^p Q^p= d^*\xi^p \ \hbox{ in } B_1,
\eeq
\beq{}\label{coulomb_gauge_p2}
d\xi^p=0  \ \hbox{ in } B_1, \  \xi^p=0 \ \hbox{ on }\partial B_1,
\eeq
and
\beq{}\label{coulomb_est}
\|\nabla Q^p\|_{L^2(B_1)}+\|\nabla\xi^p\|_{L^2(B_1)}\le C\|\Theta^p\|_{L^2(B_1)}\le C\epsilon_0.
\eeq
Multiplying  the equation (\ref{approx_dirac_hm}) by $\left(Q^p\right)^{-1}$ (see also (\ref{dirac-hm5})), we obtain
\beq{} \label{approx_dirac_hm1}
d^*\left[\left(Q^p\right)^{-1}\left(\begin{matrix}\langle d\phi_p,e_1^p\rangle \\ \vdots \\
\langle d\phi_p, e_k^p\rangle\end{matrix}\right)\right]
=-d^*\xi^p\cdot
\left(Q^p\right)^{-1}\left(\begin{matrix}\langle d\phi_p,e_1^p\rangle
\\ \vdots \\ \langle d\phi_p, e_k^p\rangle\end{matrix}\right)
+\left(Q^p\right)^{-1}\left(\begin{matrix}u_p^1\\ \vdots \\ u_p^k \end{matrix}\right).
\eeq
Since $e_i^p=\hat e_i(\phi_p)$, it is easy to see that for $1\le i\le k$,
$e_i^p\rightharpoonup e_i=\hat e_i(\phi)$ in $H^1(B_1)$  and hence $\{e_i\}$ is
an orthonormal frame along the map $\phi$.

After passing to possible subsequences, we may now assume that
$$Q^p\rightarrow Q \hbox{ weakly\ in } H^1(B_1, {\rm{SO}}(k)),
\ {\rm{ strongly\ in }}\ L^2(B_1, {\rm{SO}}(k)),  \ {\rm{ and \ a.e.\ in\ }} B_1,$$
$$ \xi^p\rightarrow \xi \hbox{ weakly\ in } H^1(B_1, {\rm{so}}(k)),
\ {\rm{ strongly\ in }}\ L^2(B_1, {\rm{so}}(k)),  \ {\rm{ and \ a.e.\ in\ }} B_1.$$
It is not hard to see that
$$\langle de_i^p, e_j^p\rangle \rightarrow \langle de_i, e_j\rangle {\rm{\ weakly\ in \ }}L^2(B_1).$$
Since $\left(\Omega_{ij}^p\right)$ is bounded in $L^2(B_1)$ and
$$\Omega_{ij}^p\rightarrow \Omega_{ij}\equiv \sum_{\alpha=1}^n \left[
\sum_{l,m=1}^k R^N(\phi)\left(e_i, e_j, e_l, e_m\right)\langle \psi^m,
\frac{\partial}{\partial x_\alpha}\circ \psi^l\rangle\right]\,dx_\alpha, \ \ {\rm{a.e.\ in \ }} B_1,$$
$\Omega_{ij}^p\rightarrow \Omega_{ij}$ weakly in $L^2(B_1)$.
Hence $\Theta_{ij}^p\rightarrow \Theta_{ij}\equiv\Omega_{ij}+\langle de_i, e_j\rangle$
weakly in $L^2(B_1)$.  Thus, sending $p\rightarrow\infty$, (\ref{coulomb_gauge_p1}) and
(\ref{coulomb_gauge_p2}) yield that $Q, \xi, \Theta$ satisfy:
\beq{}\label{coulomb_gauge_3}
Q^{-1}dQ-Q^{-1} \Theta Q= d^*\xi \ \hbox{ in } B_1,
\eeq
\beq{}\label{coulomb_gauge_4}
d\xi=0  \ \hbox{ in } B_1, \  \xi=0 \ \hbox{ on }\partial B_1.
\eeq
Since $u_p\rightarrow 0$ in $H^{-1}(B_1)$, we have that
\beq{}\label{rhs2_conv}
\left(Q^p\right)^{-1}\left(\begin{matrix}u_p^1\\ \vdots \\ u_p^k \end{matrix}\right)\rightarrow 0
\eeq
in the sense of distribution on $B_1$.
It is also easy to see
\beq{}\label{lhs_conv}
d^*\left[\left(Q^p\right)^{-1}\left(\begin{matrix}\langle d\phi_p,e_1^p\rangle \\ \vdots \\
\langle d\phi_p, e_k^p\rangle\end{matrix}\right)\right]
\rightarrow d^*\left[Q^{-1}\left(\begin{matrix}\langle d\phi,e_1\rangle \\ \vdots \\
\langle d\phi, e_k\rangle\end{matrix}\right)\right]
\eeq
in the sense of distribution on $B_1$.

Now we want to discuss the convergence of
$$A_p:=d^*\xi^p\cdot
\left(Q^p\right)^{-1}\left(\begin{matrix}\langle d\phi_p,e_1^p\rangle
\\ \vdots \\ \langle d\phi_p, e_k^p\rangle\end{matrix}\right).$$
Note that the $i$-th component of $A_p$, $A_p^i$, is given by
$$A_p^i=d^*\xi_{il}^p\cdot\left(Q^p_{ml}\langle d\phi_p, e_m^p\rangle\right)
=\left\langle d^*\xi_{il}^p\cdot d\phi_p,  \left(Q^p_{ml}e_m^p\right)\right\rangle.$$

For this, we recall a compensated compactness lemma, see Freire-M\"uller-Struwe \cite{FMS} and Wang \cite{wang1}
Lemma 3.4 for a proof.
\ble{} \label{compen-compact1} For $n=2$, suppose that $f_p\rightarrow f$ weakly in $H^1(B_1)$, $g_p\rightarrow g$
weakly in $H^1(B_1, \wedge^2\mathbb R^2)$, and $h_p\rightarrow h$ weakly in $H^1(B_1)$.
Then, after passing to possible subsequences, we have
\beq{}\label{compen-compact2}
df_p\cdot d^*g_p \cdot h_p \rightarrow df\cdot d^*g \cdot h +\nu
\eeq
in the sense of distributions on $B_1$, where $nu$ is a signed Radon measure
given by
$$\nu=\sum_{j\in J}a_j \delta_{x_j},$$
where $J$ is at most countable, $a_j\in\mathbb R$, $x_j\in B_1$, and
$\sum_{j\in J}|a_j|<+\infty$.
\ele

Applying Lemma \ref{compen-compact1}, we conclude that for $1\le i\le k$,
\beq{} \label{rhs1_conv}
A_p^i\rightarrow A^i:= d^*\xi_{il}\cdot\left(Q_{ml}\langle d\phi, e_m\rangle\right)+\nu^i \ {\rm{\ in\ } } B_1
\eeq
where
$$\nu^i=\sum_{j=1}^\infty a_j^i \delta_{x_j^i}, \ \sum_{j=1}^\infty |a_j^i|<+\infty.$$
Putting (\ref{lhs_conv}), (\ref{rhs1_conv}), and (\ref{rhs2_conv}) into (\ref{approx_dirac_hm1}), we obtain
\beq{}\label{approx_dirac_hm2}
d^*\left[Q^{-1}\left(\begin{matrix}\langle d\phi,e_1\rangle \\ \vdots \\
\langle d\phi, e_k\rangle\end{matrix}\right)\right]
=-d^*\xi\cdot
Q^{-1}\left(\begin{matrix}\langle d\phi,e_1\rangle
\\ \vdots \\ \langle d\phi, e_k\rangle\end{matrix}\right)
+\left(\begin{matrix}\nu^1\\ \vdots \\ \nu^k\end{matrix}\right).
\eeq
Note that (\ref{approx_dirac_hm2}) implies
$$\left(\begin{matrix}\nu^1\\ \vdots \\ \nu^k\end{matrix}\right)\in H^{-1}(B_1)+L^1(B_1)$$
so that $\nu^i=0$ for all $1\le i\le k$. Therefore, we have
\beq{}\label{approx_dirac_hm3}
d^*\left[Q^{-1}\left(\begin{matrix}\langle d\phi,e_1\rangle \\ \vdots \\
\langle d\phi, e_k\rangle\end{matrix}\right)\right]
=-d^*\xi\cdot
Q^{-1}\left(\begin{matrix}\langle d\phi,e_1\rangle
\\ \vdots \\ \langle d\phi, e_k\rangle\end{matrix}\right).
\eeq
This and (\ref{coulomb_gauge_3}) imply that
$$d^*(\langle d\phi, e_i\rangle)=\sum_{j=1}^k\Theta_{ij}\cdot \langle d\phi, e_j\rangle.$$
Note that this equation is equivalent to the Dirac-harmonic map equation $\tau(\phi)=\mathcal R^N(\phi,\psi)$.

Now we want to show $D\hspace{-0.25cm}/\hspace{0.1cm}\psi=0$. To see this, observe that if we write
$\psi_p=\psi_p^i\tensor e_i^p$ and $v_p=v_p^i\tensor e_i^p$, then
$D\hspace{-0.25cm}/\hspace{0.1cm}\psi_p=v_p$ becomes
\beq{}\partial\hspace{-0.2cm}/\hspace{0.1cm} \psi_p^i=-\Gamma_{jl}^i(\phi_p)\langle\frac{\partial\phi_p}{\partial x_\alpha}, e_j^p\rangle
\frac{\partial}{\partial x_\alpha}\circ \psi_p^l +v_p^i. \label{approx_dirac_eqn}
\eeq
It is easy to see that, after taking $p$ to $\infty$, (\ref{approx_dirac_eqn}) yields
$$\partial\hspace{-0.2cm}/\hspace{0.1cm}\psi^i=\Gamma_{jl}^i(\phi)\langle\frac{\partial\phi}{\partial x_\alpha}, e_j\rangle
\frac{\partial}{\partial x_\alpha}\circ \psi^l,$$
this is equivalent to $D\hspace{-0.25cm}/\hspace{0.1cm}\psi=0$. Thus the proof is complete. \qed \\

\noindent{\bf Proof of Theorem \ref{conv_dirac_hm}}.\\

Define the possible concentration set
$$C=\bigcap_{R>0}\left\{x\in M: \ \liminf_{p\rightarrow\infty}\int_{B_R(x)}\left(|\nabla\phi_p|^2+|\psi_p|^4\right)
>\epsilon_1^2\right\}.$$
Then, by a simple covering argument, we have that $C$ is at most a finite subset in $M$. By the definition, we know
that for any $x_0\in M\setminus C$, there exists $r_0>0$ and a subsequence of $(\phi_p,\psi_p)$, denoted as itself,
such that
$$\lim_{p\rightarrow\infty}\int_{B_{r_0}(x_0)}\left(|\nabla\phi_p|^2+|\psi_p|^4\right)\le\epsilon_1^2.$$
Applying Lemma \ref{small_convergence}, we conclude that $(\phi,\psi)$ is a weakly Dirac-harmonic map on
$B_{r_0}(x_0)$. Since $x_0\in M\setminus C$ is arbitrary, this implies that $(\phi,\psi)$ is a weakly
Dirac-harmonic map on $M\setminus C$. Since $C$ is at most finite, one can easily show that
$(\phi,\psi)$ is also a weakly Dirac-harmonic map on $M$. Hence Theorem \ref{2d_reg} also implies
that $(\phi,\psi)$ is a smooth Dirac-harmonic map on $M$. The proof is now complete. \qed

\bigskip
\noindent{\bf ACKNOWLEDGEMENT}

\smallskip
This work was carried out when the second author visited Department of Mathematics,
University of Kentucky under a scholarship from Shanghai Jiaotong University, P.R. China.
He would like to express his gratitude to the department for both its hospitality and
excellent research environment. The first author is partially supported by NSF 0601182.
He would like to thank Dr. L. Zhao  from Beijing Normal University for explaining 
his work \cite{zhao} to him in the summer 2008.

\end{document}